% Symbundle.tex : Version prepublication 08 octobre 2007. Corrections jan 09.
%VERSION OF july 27, 1999
\input amssym.def
\input amssym.tex
% Achtung: Bei Verwendung von jltmac werden die Fontfamilien masm und
% msbm mehrfach aufgerufen. Das kann zu Kapazitaetsschwierigkeiten
% fuehren. KH 27.2.95

\def\item#1{\vskip1.3pt\hang\textindent {\rm #1}}% THIS REPLACES KNUTH'S DEF'N

                                  % THIS REPLACES KNUTH'S DEF'N

\tolerance=300
\pretolerance=200
\hfuzz=1pt
\vfuzz=1pt

% Print out with -x300 -y50

% bei unseren Artikeln uebliches Format:
\magnification=\magstep0

%\magnification=1100 
% gibt groesseres Schriftbild

% Offsetwerte fuer Ausdrucke in Erlangen, hoffset um 0.6 in groesser machen
\hoffset=0.6in
\voffset=0.8in

%\baselineskip                 
%\hsize=5.8 true in 
\hsize=5.8 true in 
\vsize=8.5 true in
%\baselineskip
\parindent=25pt
\mathsurround=1pt
\parskip=1pt plus .25pt minus .25pt
\normallineskiplimit=.99pt

\countdef\revised=100
\mathchardef\emptyset="001F % THIS REPLACES KNUTH'S DEFINITION
\chardef\ss="19
\def\3{\ss}
\def\anf{$\lower1.2ex\hbox{"}$}
\def\frac#1#2{{#1 \over #2}}
\def\>{>\!\!>}
\def\<{<\!\!<}

\def\ssarr{\hbox to 30pt{\rightarrowfill}}
\def\sarr{\hbox to 40pt{\rightarrowfill}}
\def\arr{\hbox to 60pt{\rightarrowfill}}

\def\larr{\hbox to 60pt{\leftarrowfill}}
\def\Arr{\hbox to 80pt{\rightarrowfill}}

{}

\def\ad{\mathop{\rm ad}\nolimits}

\def\Aut{\mathop{\rm Aut}\nolimits}

\def\D{{\rm D}}
\def\Der{\mathop{\rm Der}\nolimits}

\def\det{\mathop{\rm det}\nolimits}

\def\End{\mathop{\rm End}\nolimits}

\def\Gl{\mathop{\rm Gl}\nolimits}

\def\Hom{\mathop{\rm Hom}\nolimits}%
\def\id{\mathop{\rm id}\nolimits} % USED FOR IDENTITY FUNCTION

% USED FOR IMAGINARY PART OF COMPLEX NUMBERS

% USED FOR REAL PART OF COMPLEX NUMBERS

%\def\reach{\mathop{\rm reach}\nolimits}
%\def\SAut{\mathop{\rm SAut}\nolimits}
%\def\SEnd{\mathop{\rm SEnd}\nolimits}

%\def\HSp{\mathop{\rm HSp}\nolimits}
\def\Sp{\mathop{\rm Sp}\nolimits}

% USED FOR TRACE OF MATRIX
% USED FOR TRACE OF MATRIX
%\def\trdeg{\mathop{\rm trdeg}\nolimits} 

%\def\Vequiv{\mathrel{_V\!\equiv}}

\def\pr{\mathop{\rm pr}\nolimits}

\def\0{{\bf 0}}
\def\1{{\bf 1}}

\def\b{{\frak b}}

\def\f{{\frak f}}
\def\g{{\frak g}}
\def\gl{{\frak {gl}}}
\def\h{{\frak h}}

\def\k{{\frak k}}

\def\m{{\frak m}}

\def\q{{\frak q}}

\def\str{{\frak {str}}}

\def\A{{\Bbb A}}
 
\def\C{{\Bbb C}} 
\def\DD{{\Bbb D}} 
 
\def\H{{\Bbb H}} 
\def\K{{\Bbb K}}

\def\R{{\Bbb R}} 
 
\def\Z{{\Bbb Z}} 

\def\:{\colon}  %8.5.92
\def\.{{\cdot}}
\def\|{\Vert}
\def\bsk{\bigskip}

\def\giantskip{\vskip2\bigskipamount}
\def\gsk{\giantskip}

\def\msk{\medskip}

\def\ssk{\smallskip}

\def\bbr{\bigbreak}
\def\giantbreak{\par \ifdim\lastskip<2\bigskipamount \removelastskip
         \penalty-400 \giantskip\fi}

\def\nin{\noindent}
\def\cen{\centerline}
\def\pagebreak{\vskip 0pt plus 0.0001fil\break}
\def\linebreak{\break}

\def\eps{\varepsilon}
\def\epsilon{\varepsilon}

\def\nin{\noindent}

\def\pder#1,#2,#3 { {\partial #1 \over \partial #2}(#3)}
\def\pde#1,#2 { {\partial #1 \over \partial #2}}
\def\phi{\varphi}

% Besser \ltimes und \rtimes aus dem AMS-Symbols verwenden. 
%\def\sdir#1{\hbox{$\mathrel\times{\hskip -4.6pt 
%            {\vrule height 4.7 pt depth .5 pt}}\hskip 2pt_{#1}$}}

\def\tilde{\widetilde}

\font\eightrm=cmr8

% SANS SERIF 10 POINT
 %SANS SERIF 10 POINT ITALIC

%\font\smc8=cmcsc8 
 %SLANTED TYPEWRITER 10 POINT
 %BOLD FACE MATH SYMBOLS 10 POINT
 %DUNHILL STYLE 10 POINT
 %SAN SERIF BOLD EXTENDED 10 POINT
 %USED FOR TITLES
 %USED FOR TITLES
\font\bfone=cmbx10 scaled\magstep1 %BOLDFACE AT MAGSTEP 1
\font\bftwo=cmbx10 scaled\magstep2 %BOLDFACE AT MAGSTEP 2
 %BOLDFACE AT MAGSTEP 3

\def\qed{{\unskip\nobreak\hfil\penalty50\hskip .001pt \hbox{}\nobreak\hfil
          \vrule height 1.2ex width 1.1ex depth -.1ex
           \parfillskip=0pt\finalhyphendemerits=0\medbreak}\rm}
%This is the end-of-proof sign. 
%Not to be used in display mode. 
%If you want to conclude a proof 
%at the end of a line is display mode use 

%BUT OMIT $$---the macro will write that

\def\Lemma #1. {\bigbreak\vskip-\parskip\noindent{\bf Lemma #1.}\quad\it}

\def\Sublemma #1. {\bigbreak\vskip-\parskip\noindent{\bf Sublemma #1.}\quad\it}

\def\Proposition #1. {\bigbreak\vskip-\parskip\noindent{\bf Proposition #1.}
\quad\it}

\def\Corollary #1. {\bigbreak\vskip-\parskip\nin{\bf Corollary #1.}
\quad\it}

\def\Theorem #1. {\bigbreak\vskip-\parskip\noindent{\bf Theorem #1.}
\quad\it}

\def\Definition #1. {\rm\bigbreak\vskip-\parskip\noindent{\bf Definition #1.}
\quad}

\def\Remark #1. {\rm\bigbreak\vskip-\parskip\noindent{\bf Remark #1.}\quad}

\def\Example #1. {\rm\bigbreak\vskip-\parskip\noindent{\bf Example #1.}\quad}

\def\Problems #1. {\bigbreak\vskip-\parskip\noindent{\bf Problems #1.}\quad}
\def\Problem #1. {\bigbreak\vskip-\parskip\noindent{\bf Problem #1.}\quad}

\def\Conjecture #1. {\bigbreak\vskip-\parskip\noindent{\bf Conjecture #1.}\quad}

\def\Proof#1.{\rm\par\ifdim\lastskip<\bigskipamount\removelastskip\fi\smallskip
            \noindent {\bf Proof.}\quad}

\def\Axiom #1. {\bigbreak\vskip-\parskip\noindent{\bf Axiom #1.}\quad\it}

\def\Satz #1. {\bigbreak\vskip-\parskip\noindent{\bf Satz #1.}\quad\it}

\def\Korollar #1. {\bbr\vskip-\parskip\nin{\bf Korollar #1.} \quad\it}

\def\Bemerkung #1. {\rm\bigbreak\vskip-\parskip\noindent{\bf Bemerkung #1.}
\quad}

\def\Beispiel #1. {\rm\bigbreak\vskip-\parskip\noindent{\bf Beispiel #1.}\quad}
\def\Aufgabe #1. {\rm\bigbreak\vskip-\parskip\noindent{\bf Aufgabe #1.}\quad}

\def\Beweis#1. {\rm\par\ifdim\lastskip<\bigskipamount\removelastskip\fi
           \smallskip\noindent {\bf Beweis.}\quad}

\nopagenumbers

\def\date{ }

\def\title{Title ??}
\def\author{Author ??}

\def\thanks#1{\footnote*{\eightrm#1}}

\def\rightheadline{\hfil{\eightrm\author}\hfil\tenbf\folio}
\def\leftheadline{\tenbf\folio
\hfil{\eightrm\title}
\hfil}
\headline={\vbox{\line{\ifodd\pageno\rightheadline\else\leftheadline\fi}}}

\def\firstheadline{}
\def\firstfootline{\cen{\rm\folio}}

\def\seite #1 {\pageno #1
               \headline={\ifnum\pageno=#1 \firstheadline
               \else\ifodd\pageno\rightheadline\else\leftheadline\fi\fi}
               \footline={\ifnum\pageno=#1 \firstfootline\else{}\fi}}

%%%THIS IS THE MACRO LEFTSPACE.TEX %%%TO THD VIA WAFRUPP
\newdimen\dimenone
 \def\checkleftspace#1#2#3#4{%DIESER MACRO STAMMT VON APPELT
 \dimenone=\pagetotal%#1=Skip vorher,#2=Font,#3=Text,#4=Skip nachher  
 \advance\dimenone by -\pageshrink   %testen ob Titel noch mit Gewalt auf Seite 
                                                                          %geht
 \ifdim\dimenone>\pagegoal          %nacha tua nix-- gewoehnliche Outputroutine 
   \else\dimenone=\pagetotal
        \advance\dimenone by \pagestretch
        \ifdim\dimenone<\pagegoal
          \dimenone=\pagetotal
          \advance\dimenone by#1         %addieren Skip vor Ueberschrift (=#1)
          \setbox0=\vbox{#2\parskip=0pt                %#2 ist gewaehlter Font
                     \hyphenpenalty=10000
                     \rightskip=0pt plus 5em
                     \noindent#3 \vskip#4}    %#3=Ueberschrift,#4=skip nachher
        \advance\dimenone by\ht0
        \advance\dimenone by 3\baselineskip   
        \ifdim\dimenone>\pagegoal\vfill\eject\fi
          \else\eject\fi\fi}

%%% OUR HEADLINE MACROS LOOK LIKE THIS USING THIS MACRO

\def\subheadline #1{\nin\bigbreak\vskip-\lastskip
      \checkleftspace{0.7cm}{\bf}{#1}{\medskipamount}
          \indent\vskip0.7cm\centerline{\bf #1}\medskip}

\def\sectionheadline #1{\bigbreak\vskip-\lastskip
      \checkleftspace{1.1cm}{\bf}{#1}{\bigskipamount}
         \vbox{\vskip1.1cm}\cen{\bfone #1}\bsk}

\def\lsectionheadline #1 #2{\bigbreak\vskip-\lastskip
      \checkleftspace{1.1cm}{\bf}{#1}{\bigskipamount}
         \vbox{\vskip1.1cm}\cen{\bfone #1}\msk \cen{\bfone #2}\bsk}

\def\lchapterheadline #1 #2{\bigbreak\vskip-\lastskip\indent\vskip3cm
                       \cen{\bftwo #1} \msk \cen{\bftwo #2} \gsk}
\def\llsectionheadline #1 #2 #3{\bigbreak\vskip-\lastskip\indent\vskip1.8cm
\cen{\bfone #1} \msk \cen{\bfone #2} \msk \cen{\bfone #3} \nobreak\bsk\nobreak}

%\def\[#1 #2\par{\hbox{\vtop{\hsize = 2.5 true cm \nin [#1]\hfill}
%\vtop{\hsize = 12.0 true cm \nin #2\penalty10000\llap.}}
%\vbox{\vskip.3\baselineskip}}  
% parameters for hsize are percentages !!! f.e. 0.2 + 0.8 = 1.0

\newtoks\literat
\def\[#1 #2\par{\literat={#2\unskip.}%
\hbox{\vtop{\hsize=.15\hsize\nin [#1]\hfill}
\vtop{\hsize=.82\hsize\nin\the\literat}}\par
\vskip.3\baselineskip}

\def\references{
\sectionheadline{\bf References}
\frenchspacing

\entries\par}

\mathchardef\emptyset="001F 
\def\address{Author: \tt$\backslash$def$\backslash$address$\{$??$\}$}

\def\firstpage{\nin
{\obeylines \parindent 0pt }
\vskip2cm
\centerline{\bfone\title}
\gsk
\centerline{\bf\author}
\vskip1.5cm \rm}

\def\lastpage{\par\vbox{\vskip1cm\nin
\line{
\vtop{\hsize=.5\hsize{\parindent=0pt\baselineskip=10pt\nin\address}}
\hfill} }}

% END OF LIEMACS.TEX

\def\DD{{\rm D}}
\def\D{{\Bbb D}}

\def\Vaut{{\rm VertAut}}

\def\date{otober 7, 2007}
\overfullrule = 0pt 

 \pageno=1

\def\title{Symmetric bundles and representations of
Lie triple systems}

% F\"ur endg\"ultige Version dieses Makro l\"oschen
%\def\leftheadline{\tenbf\folio\hfil Symbundle.tex, \eightrm
%prelimary version, \date}

\def\author{Wolfgang Bertram, Manon Didry }

\def\address
{}

\firstpage

\centerline
{Institut Elie Cartan -- U.M.R. 7502}
\centerline
{Universit\'e Henri Poincar\'e (Nancy I)}
\centerline
{B.P. 239}
\centerline
{54506 Vand\oe uvre-l\`es-Nancy Cedex}
\centerline
{France}

\ssk
\centerline
{e-mail: {\tt bertram@iecn.u-nancy.fr}, $\,$
{\tt didrym@iecn.u-nancy.fr}}

\bigskip
 
\nin {\bf Abstract.} We define
{\it symmetric bundles} as vector bundles in the category of symmetric spaces;
it is shown that this notion is the geometric analog of the one of a {\it 
representation of a Lie triple system}.
A symmetric bundle has an {\it underlying reflection space}, and 
we investigate the corresponding forgetful functor both from the point of
view of differential geometry and from the point of view of representation theory.
This functor is not injective, as is seen 
by constructing ``unusual'' symmetric bundle structures on the tangent bundles
of certain symmetric spaces. 

\msk \nin {\bf AMS subject classification:}
17A01, % Non assoc alg: General theory
17B10, % Lie algebras:Representations, algebraic theory (weights)
% 53B05, % local diff geo : Linear and affine connections
53C35 % Symmetric spaces 

\msk \nin {\bf Keywords:} symmetric space, symmetric bundle,
general representations, Lie triple system, canonical connection

%\vskip 20mm

\sectionheadline{Introduction}

Although this is not common, {\it linear representations of Lie groups}
may be defined as {\it vector bundles in the category of Lie groups}:
 if $\rho:G \to \Gl(V)$ is a (say, finite-dimensional) representation
of a Lie group in the usual sense, then the semidirect product
$F:=G \ltimes  V$ of $G$ and $V$ is a Lie group and at
the same time a vector bundle over $G$ such that both structures are compatible in
the following sense:

\ssk
\item{(R1)}  the projection $\pi:F \to G$ is a Lie group homomorphism,
\item{(R2)}  the group law $F \times F \to F$ is a morphism of vector bundles,
i.e., it preserves fibers and, fiberwise, group multiplication
 $F_g \oplus F_h \to F_{gh}$, $(v,w) \mapsto vw$ is linear.

\ssk \nin Conversely,  given a vector bundle $F$ over $G$ 
with total space a Lie group and
having such properties,
the representation of $G$ can be recovered as the fiber $F_e$  over the unit element
$e$ on which $G$ acts by conjugation. 
For instance, the tangent bundle $TG$ corresponds to the adjoint representation, and the
cotangent bundle $T^*G$ to the coadjoint representation of $G$.

\msk
In this work we wish to promote the idea that this way of viewing representations
is the good point of view when looking for a notion of ``representation'' for
other categories of spaces which, like Lie groups, are defined by one or several
``multiplication maps'': somewhat simplified, {\it a representation of a given
object $M$ of such a category is a vector bundle over $M$ in the given category},
where ``vector bundle in the given category'' essentially means that the analogues of
(R1) and (R2) hold. In fact, this simple notion came out as a result of our
attemps to find a ``global'' or ``geometric'' analog of the notion of
{\it representation of general $n$-ary algebraic structures}:
given a multilinear algebraic structure defined by identities, such as Lie-, Jordan-
or other algebras or triple systems,  Eilenberg [Ei48] introduced  a natural
notion of
{\it (general) representation} (which is widely used in Jordan theory, see
[Jac51], [Lo73], [Lo75]).\footnote{$^1$}{\eightrm
A word of warning: in the literature, especially on Jordan algebras, there is
some confusion in terminology; the notion of general representation differs very much
from the idea of a representation to be a homomorphism ``into some matrix realization''.
Unfortunately, the word ``representation'' is also used in this second sense for
Jordan algebras (cf., e.g., [FK94]) and for symmetric spaces ([Be00, I.5]); we suggest
to replace this by the term ``specialization'', in the sense of ``homomorphism into
a special (i.e., matrix or operator) object''.}
Essentially, a representation $V$ of such an $n$-ary algebra $\m$ is equivalent to
defining on the direct sum $V \oplus \m$ an $n$-ary algebraic structure satisfying the
same defining identities as $\m$ and such that some natural properties hold, which
turn out to be exactly the ``infinitesimal analogs'' of (R1) and (R2): for instance,
$V$ will be an ``abelian'' ideal in $V \oplus \m$, corresponding to the role of the fiber in
a vector bundle. The archetypical example is given by the {\it adjoint representation}
which is simply $\m \oplus \epsilon \m$ with $\epsilon^2=0$, the scalar extension of
$\m$ by dual numbers, which shall of course correspond to the tangent bundle in the
geometric picture. However, nothing guarantees in principle that there be a
``coadjoint representation'' and a ``cotangent bundle in the given category''!

\msk
This approach is very general and has a  wide range of possible applications:
at least locally, any affine connection on a manifold gives rise to a smooth
``multiplication map'' (a {\it local loop}, see [Sab99]), which by deriving gives rise to
$n$-ary algebras, and hence may be ``represented'' by vector bundles.   
Concretely, we will show how all these ideas work for the most proeminent example
 of such structures, namely for {\it symmetric spaces}
(here, the approach to symmetric spaces by Loos [Lo69] turns out to be best suited;
we recall some basic facts and the relation with homogeneous spaces $M=G/H$ in Chapter 1).
{\it Symmetric bundles} are defined as vector bundles in the category of symmetric
spaces (Section 1.4) and
their infinitesimal analogs, {\it representations of
Lie triple systems} are introduced (Chapter 2). These have already been studied from
a purely algebraic point of view by  T.\ Hodge and B.\ Passhall [HP02].
Another algebraic point of view ([Ha61]) features the aspect of {\it representations of 
Lie algebras with involution} (cf.\ Section 3.2),  which we use to prove that, in the real 
finite-dimensional case, such representations are in one-to-one 
correspondence with symmetric bundles (Theorem 3.4).
This result implies that the cotangent bundle $T^*M$ of a real finite-dimensional symmetric
space is again a symmetric bundle -- this is much less obvious than the corresponding
fact for the tangent bundle, and it owes its validity to the fact that, on the
level of representations of Lie triples systems, every representation admits a
{\it dual representation} (Section 4.2).
Whereas, among the algebraic constructions of new vector bundles from old ones,
the {\it dual} and the {\it direct sum constructions}
 survive  in the category of symmetric bundles, 
this is not the case for
{\it tensor products} and {\it hom-bundles}: they have to replaced by other, more
complicated constructions (Chapter 4).

\msk
Compared to the case of
 group representations, a new feature of symmetric bundles
is that they are ``composed ojects'': for a Lie group, the group structure on
$G \ltimes V$ and its structure of a {\it homogeneous vector bundle} 
are entirely equivalent.  For symmetric bundles, the structure of a homogeneous
vector bundle carries strictly less information than that of the symmetric bundle:
let us assume $F$ is a symmetric bundle over a {\it homogeneous} symmetric space
$M=G/H$; then $F$ carries two structures: it is a homogeneous symmetric space $F=L/K$,
and, under the action of the smaller group $G$, it is
 a homogeneous vector bundle $G \times_H V$, with
$V=F_o$ being the fiber over the base point $o=eH$.
Basically, seeing $F$ as a
homogeneous vector bundle only retains the representation of $H$ on $F_o$,
whereas seeing $F$ as symmetric space $L/K$ takes into account 
 the whole isotropy representation 
of the bigger group $K$. In other words, there is a 
{\it forgetful functor from symmetric bundles to homogeneous vector bundles}.
Conversely, the following ``extension problem'' arises: {\it given a homogeneous
symmetric space $M=G/H$,
which homogeneous vector bundles (i.e., which $H$-representations) admit a
compatible structure of a symmetric bundle ?}
%
%More conceptually and more generally, our forgetful functor is interpreted as a
%{\it forgetful functor from symmetric bundles to reflection spaces} (Proposition 1.7).
%
On the infinitesimal level of general
representations of Lie triple systems, the forgetful functor appears as follows:
a general representation of a Lie triple system $\m$ consists of {\it two}
trilinear maps $(r,m)$, and we simply forget
the second component $m$ (Section 2.9). 
The extension problem is then: {\it when does $r$ admit a compatible trilinear map $m$
such that $(r,m)$ is a representation of Lie triple systems?}
For a geometric interpretation of this problem, one notes that
the trilinear map $r$ is of the type of a {\it curvature tensor}, and indeed 
one can prove that every symmetric bundle admits a canonical connection (Theorem 5.1)
such that $r$ becomes its curvature tensor (Theorem 5.3). 
It seems thus that the representations of $H$ that admit an extension to a symmetric
bundle are those that can themselves be interpreted as holonomy
representation of a connection on a vector bundle. However, the situation is
complicated by the fact that the canonical connection on $F$ does {\it not} 
determine completely the symmetric bundle structure on $F$.
% (but a derived connection on $TF$ does). 

\msk
We do not attack in this work the problem of classifying representations of, say,
finite dimensional simple symmetric spaces; but we give a large class of examples
of ``unusual'' symmetric bundle structures on tangent bundles (Chapter 6), thus showing
that  the above mentioned forgetful functor is not injective.
 In fact, as observed in [Be00], many (but not all)
 symmetric spaces $M=G/H$ admit, besides their ``usual'' complexification
$M_\C =G_\C/H_\C$, another, ``twisted'' or ``hermitian'' one $M_{h\C}=L/K$.
We show here that a similar construction works when one replaces
``complexification'' by ``scalar extension by dual numbers'' (replace the condition
$i^2=-1$ by $\eps^2 = 0$), and that in this way we obtain two different symmetric
structures on the tangent bundle $TM$. 
A particularly pleasant example is the case of the general linear group
$M=\Gl(n,\R)$: in this case, the usual tangent bundle $TM$ is the group
 $\Gl(n,\R[\epsilon])$
(scalar extension by dual numbers), whereas the ``unusual'' symmetric structure on
the tangent bundle is obtained by realizing $TM$ as the homogeneous space 
$\Gl(n,\D)/\Gl(n,\R[\eps])$, where $\D$ is some
degenerate version of the quaternions (Theorem 6.2). 
We conjecture that, for real simple symmetric spaces, there are no other
symmetric bundle structures on the tangent space than the ones just mentioned. 
In other words, we conjecture
that the extension problem as formulated here is closely related to the
``extension problem for the Jordan-Lie functor'' 
from [Be00]; however, this remains a topic for future research.

\msk
The results presented in this paper partially extend results from the 
thesis [Did06], where a slightly different axiomatic definition of symmetric bundles
was proposed in a purely algebraic setting, permitting to state the analog of
Theorem 3.4 (equivalence of symmetric bundles and representations of Lie triple systems)
in an algebraic framework (arbitrary dimension and arbitrary base field; see 
Theorem 2.1.2 in loc.\ cit.), based on results published in [Did07].
The present paper is independent from the results of [Did07], but nevertheless the
framework still is quite general: our symmetric spaces are of arbitrary dimension and
defined over very general topological base fields or rings $\K$ -- 
for instance, the setting includes real or complex
infinite dimensional (say, Banach)
 symmetric spaces or $p$-adic symmetric spaces (Section 1.1).
We hope the reader will agree that, in the present case, this degree of generality
does not complicate the theory, but rather simplifies it by forcing one to
search for the very basic concepts.

\ssk
After this work had been finished, we learned from Michael Kinyon that
 the question of defining ``modules'' for an
object in a category had already been investigated by J.\ M.\ Beck in his thesis
([Beck67]; see also [Barr96]): he defines a {\it module} to be an abelian group object
in the slice category over the given object. It seems reasonable to conjecture
that, in the cases considered here, this notion should agree with ours, but by
lack of competence in category theory we have not been able to check this.

\sectionheadline
{1. Symmetric bundles}

\nin {\bf 1.1. Notation and general framework.}
This work can be read on two different levels:
the reader may take $\K=\R$ to be the real base field and understand 
by ``manifold'' finite-dimensional real manifolds in the usual sense;
then our symmetric spaces and Lie groups are the same as in
[Lo69] or [KoNo69], or one may consider a commutative topological
field or ring $\K$, having dense unit group $\K^\times$ and such
that $2$ is invertible in $\K$; then we refer to [Be08] for the
definition of manifolds and Lie groups  over $\K$.
Readers interested in the general case should just
keep in mind that, in general,

\ssk
\item{--} symmetric spaces need no longer be homogeneous (cf.\ Item 1.3 (2) below),
\item{--} there is no exponential map and hence no general tool to ``integrate'' 
infinitesimal structures to local ones.

\nin
\nin If we use such tools, it will be specifically mentioned that we are in the
real (or complex) finite dimensional case. --
In the sequel, the word {\it linear space} means ``(topological) $\K$-module''.

\msk
\nin {\bf 1.2. Symmetric spaces and reflection spaces.}
A {\it reflection space} (``Spiegelungsraum'', introduced by O.\ Loos in [Lo67])
 is a smooth manifold $M$ together with a smooth
``product map'' $\mu: M \times M \to M$, $(x,y) \mapsto \mu(x,y)=:
\sigma_x(y)$ satisfying, for all $x,y,z \in M$,

\msk

\item{(S1)}
$\mu(x,x)=x$,
\item{(S2)}
$\mu(x,\mu(x,y))=y$,
\item{(S3)}
$\sigma_x \in \Aut(\mu)$, i.e. $\mu(x, \mu(y,z))=\mu(\mu(x,y),\mu(x,z))$.

\msk
\nin
The reflection space $(M,\mu)$ 
is called a {\it symmetric space} if in addition

\item{(S4)}
for all $x \in M$, the differential $T_x(\sigma_x)$ of the ``symmetry'' $\sigma_x$ at
$x$ is the negative of the identity of the tangent space $T_x M$.

\msk \nin In the real finite-dimensional case this is (via the implicit function
theorem) equivalent to

\item{(S4')}
for all $x \in M$, the fixed point $x$ of $\sigma_x$ is isolated.

\ssk
\nin
{\it Homomorphisms} of such structures are smooth maps which commute
with  product maps. According to (S3), all maps of the form
$\sigma_x \circ \sigma_y$, $x,y \in M$, are automorphisms; the subgroup
$G(M)$ of $\Aut(M)$ generated by these elements is called the {\it transvection
group of $M$.}
 Often one considers the category of reflection
spaces, resp.\ symmetric spaces {\it with base point}:
a {\it base point} is just a distinguished point, often denoted by
$x_0$ or $o$, and homomorphisms are then required to preserve base points.
If $o\in M$ is a base point, one defines the {\it quadratic map}
by
$$
Q:=Q_o:M \to \Aut(M), \quad x \mapsto Q(x):=\sigma_x \circ \sigma_o
$$
and the {\it powers} by 
$x^{2k}:=Q(x)^k(o)$  and  $x^{2k+1}:=Q(x)^kx$.

\msk
\nin {\bf 1.3. Examples.}
(1) The group case. Every Lie group with the new multiplication $\mu(g,h)=gh^{-1}g$
is a symmetric space.

\ssk
(2) Homogeneous symmetric spaces. We say that a symmetric space is {\it homogeneous}
if the group $G:=G(M)$ acts transitively on it and carries a Lie group structure such
that this action is smooth. Let $o$ be a base point and $H$ its stabilizer, so that
$M\cong G/H$. Then the map $\sigma:G \to G$, $g \mapsto \sigma_o \circ g \circ 
\sigma_o$ is an involution of $G$, and the multiplication map on $G/H$ is given by
$$
\sigma_o(gH)=\sigma(g)H, \quad \quad \sigma_{gH}(g'H)=g \sigma(g)^{-1} g'H.
$$
In finite dimension over $\K=\R$, every connected symmetric space is of this form,
for a suitable involution $\sigma$ of a Lie group $G$
(see [Lo67], [Lo69]).

\ssk
(3) Linear symmetric spaces.
Assume $V$ is a linear space over $\K$; we consider $V \times V$
as a linear space and thus write $V\oplus V$.
Assume that $V$ carries a symmetric 
space structure $\mu:V \oplus V \to V$ which is
a $\K$-linear map.
Because of (S4), the symmetry $s_0 = \mu(0,\cdot):V \to V$, being
a linear map, must agree with its tangent map $-\id_V$. Then it follows
that
$$
\eqalign{
\mu(v,w) & =\mu((v,v)-(0,v)+(0,w)) = \mu(v,v) - \mu(0,v) + \mu(0,w) \cr
&  =
v - (-v) - w = 2v-w. \cr}
$$
Conversely, every linear space equipped with the multiplication map $\mu(v,w)=2v-w$
is a symmetric space. (In fact, it is the group case $G=V$.)
 With respect to the zero vector as base point,
$Q(x)=\tau_{2x}$ is translation by $2x$, and the powers are
$x^k=kx$. 

\ssk
(4) Polynomial symmetric spaces. 
In the same way as in the preceding example, we can consider linear spaces together
with a symmetric structure which is a {\it polynomial} map $V \oplus V \to V$
 -- see [Did06] for a theory of such spaces.

\msk
\nin {\bf 1.4. Symmetric  bundles.}
A {\it symmetric bundle} (or, longer but more precise: {\it symmetric vector bundle})
 is a vector bundle $\pi:F \to M$ such that
\ssk

\item{(SB1)}
$(F,\mu)$ and $(M,\mu_M)$ are symmetric spaces such that
$\pi:F \to M$ is a homomorphism of symmetric spaces,
\item{(SB2)}
for all $(p,q) \in M \times M$, the map induced by $\mu:F \times F \to F$ fiberwise,
$$
F_q \oplus F_p \to F_{\mu(p,q)}, \quad (v,w) \mapsto \mu(v,w)
$$
(which is well-defined according to (SB1)),
is linear. 

\ssk \nin
{\it Homomorphisms} of symmetric 
bundles are vector bundle homomorphisms
that are also homomorphisms of the symmetric spaces in question.
Clearly, the concept of symmetric bundle could be adapted to other classes
of bundles whenever the fibers belong to a category that admits direct products
(e.g., {\it multilinear bundles} in the sense of [Be08]):
it suffices to replace (SB2) by the requirement
 that the map $F_q \times F_p \to F_{\mu(p,q)}$ be a morphism in that category.
Also, it is clear that such concepts exist for any category of manifolds equipped
with binary, ternary or other ``multiplication maps'', such as {\it generalized
projective geometries} (cf.\ [Be02]).
 However, in the sequel we will stick to the
case of vector bundles and symmetric spaces.
%
%[see, however, some remarks on the Jordan case, at the end?]

\ssk
A symmetric bundle is called {\it trivial} if it is trivial as a bundle, and if,
as a symmetric space, it is simply
 the direct product of $M$ with a vector space.
The first non-trivial example of a symmetric  bundle is the tangent bundle
$F:=TM$ of a symmetric space $(M,\mu)$: as to (SB1), it is well-known that
$TM$ with product map $T\mu:TM \times TM \to TM$ is a symmetric space
such that the canonical projection is a homomorphism and the fibers
are flat subspaces
(see [Lo69] for the real finite dimensional and [Be08] for the general
case). Property (SB2) follows immediately
 from the linearity of the tangent map
$T_{(p,q)}\mu:T_p M \times T_q M \to T_{\mu(p,q)}M$.

\msk \nin
{\bf 1.5. Some elementary properties of symmetric bundles.} For a symmetric bundle
$F$ over $M$, the following holds:

\ssk
\item{(SB3)}
the symmetric space structure on the fiber $F_x=\pi^{-1}(x)$ over
$x \in M$ coincides
with the canonical symmetric space structure of the vector space
$(F_x,+)$ (i.e., $\mu(u,v)=2u-v$).

\item{(SB4)} The zero-section $z:M \to F$ is a homomorphism of symmetric
spaces. (Hence in the sequel we may identify $M$ with $z(M)$, and the use of the
same letter $\mu$ for the multiplication maps of $M$ and $F$ does not lead to
confusion.)
\item{(SB5)} For all $r \in \K$,   the fiberwise dilation map 
$$
(r)_F: F \to F, \quad v \mapsto rv
$$
is an endomorphism of the symmetric space $F$; for $r \in \K^\times$
it is an automorphism.
\ssk

\nin In fact, 
for $p=q$, (SB2) says that the fiber $F_p$ is a symmetric subspace
of $F$ such that its structure map $F_p \oplus F_p \to F_p$ is linear, and (SB3) 
now follows in view of Example 1.3 (3).
Since a linear map sends zero vector to zero vector,
we have
$$
\mu(0_p,0_q)=0_{\mu(p,q)},
$$
proving (SB4), and to prove (SB5), just note that 
for $v \in F_p$ and $w \in F_q$, by (SB2),
$$
(r)_F \big(\mu(v,w)\big) =r \mu(v,w) = 
\mu (rv,rw) = \mu \big( (r)_F v,(r)_F w \big).
$$
In particular, note that
 $(0)_F$ is the projection onto the zero section, and that
$(-1)_F$ can be seen as a ``horizontal reflection with respect to the zero section''.

\msk \nin {\bf 1.6. Horizontal and vertical symmetries.}
Let $\pi:F \to M$  a symmetric bundle and $u \in F_p$.
We define the {\it horizontral} (resp., {\it vertical}) 
{\it symmetry (with respect to $u$)} by
$$
\eqalign{
\vartheta_u & := \sigma_{u \over 2} \circ (-1)_F \circ \sigma_{u\over 2} \cr
\nu_u &:= \sigma_{u \over 2} \circ \sigma_{0_p} \circ (-1)_F \circ \sigma_{u\over 2}.
\cr}
$$
For $u=0_p$, the maps $\sigma_u,\vartheta_u$ commute with each
other because of (SB5):
$$
(-1)_F \circ \sigma_{0_p} \circ (-1)_F = \sigma_{-0_p} = \sigma_{0_p}.
$$
Therefore $\nu_u$ then is also is of order $2$.
Conjugating by $\sigma_{u \over 2}$, we see that for all $u \in F_p$, 
we get three pairwise commuting automorpisms $\sigma_u, \vartheta_u, \nu_u$,
 of order $2$ and fixing the point $\sigma_{u \over 2}(0_p)=u$.

\Lemma 1.7.
 The vertical symmetry depends only on the fiber $F_p$, that is,
for all $u,w \in F_p$, we have $\nu_u=\nu_w = \nu_{0_p} = (-1)_F \circ \sigma_{0_p}$.

\Proof.
Let us show that $\nu_u = \nu_0$ with $0=0_p$, i.e.,
%
%$$\sigma_{u \over 2} \circ \nu_0 \circ \sigma_{u\over 2} = \nu_0$$
%or
$$
\sigma_{u \over 2} \circ (-1)_F\circ \sigma_0 \circ \sigma_{u\over 2} = 
(-1)_F \circ \sigma_0.
$$
But this follows from $\sigma_0 \sigma_{u\over 2}\sigma_0 = 
\sigma_{\sigma_0({u \over 2})} =
\sigma_{-{u\over 2}}$ and the fact that $(-1)_F$ is an automorphism.
\qed

\Proposition 1.8. 
The space $F$ together with the binary map
$\nu:F \times F \to F$, $(v,w) \mapsto \nu_v(w)$ is a reflection space.

\Proof.
The defining
properties (S1) and (S2) say that $\nu_v$ is of order $2$ and fixes $v$, and
this has already been proved above. In order to establish (S3),
let $v \in F_p$ and $w \in F_q$. Then, using the preceding lemma,
$$
\eqalign{
\nu_v \nu_w \nu_v & = \nu_{0_p}\nu_{0_q}\nu_{0_p} =
(-1)_F \circ \sigma_{0_p} \circ (-1)_F \circ \sigma_{0_q} \circ (-1)_F
\circ  \sigma_{0_p} \cr
& =
(-1)_F \circ \sigma_{0_p} \circ \sigma_{0_q} \circ\sigma_{0_p} =
(-1)_F \circ  \sigma_{\mu(p,q)} \cr
& =
\nu_{\nu_v(w)} .\cr}
$$
\qed

\nin We say that $(F,\nu)$ is the {\it reflection space associated to
the symmetric bundle $(F,\mu)$}. 
Thus we get a functor $(F,\mu) \mapsto (F,\nu)$ from symmetric vector bundles
to reflection spaces; it will be a recurrent theme in this work to interprete
this functor as a forgetful functor.
Note that the differential of
$\nu_u$ has $1$-eigenspace tangent to the fiber through $u$ and
$-1$-eigenspace complementary to it; thus the distribution of the ``vertical''
 $1$-eigenspaces
is integrable, whereas the distribution of the ``horizontal''
 $-1$-eigenspaces is in general not
(see Chapter 5: the curvature of the corresponding Ehresmann connection does 
in general not vanish).

\msk \nin {\bf 1.9. Automorphisms downstairs and upstairs.}
The canonical projection $\pi:\Aut(F) \to \Aut(M)$, $(\tilde g,g) \mapsto g$
does in general not admit a cross-section; we cannot even guarantee that it is
surjective. However, it is easily seen that the projection
of transvection groups $G(F) \mapsto G(M)$ is surjective:  
 write $g \in G(M)$ as a composition of symmetries at points of
$M$; identifying $M$ with the zero section in $F$ we see that
$g$ gives rise to an element $\tilde g \in G(F)$ with $\pi(\tilde g,g)=g$.
In particular, if $M$ is homogeneous, then so is $F$:
in fact, if $v \in F_x$, then there exists $g \in G(M)$ 
with $g.o=x$; then  $\sigma_{v\over 2} \circ \tilde g(0_o)=v$.
In the real finite-dimensional case, we may replace $G(M)$ by its universal
covering; then the zero section $z:M \to F$ induces a homomorphism 
of this universal covering into $G(F)$, having discrete kernel. 
Hence, if we write $F=L/K$,
 it is not misleading to think of $G$ as a subgroup of $L$ and of $H$ as
a subgroup of $K$ (possibly up to a discrete subgroup).

\msk \nin {\bf 1.10. Homogeneous bundles over symmetric spaces.}
Assume that $M=G/H$ is a homogeneous symmetric space (Example 1.3 (2)).
To any smooth action $H \times U \to U$ on a manifold $U$ one can associate the
homogeneous bundle
$$
F= G \times_H U = G \times U/\sim, \quad \quad
(gh,v) \sim (g,hv) \, \, \forall h \in H.
$$ 
When the base $M=G/H$ is a symmetric space, we define
$\nu:F \times F \to F$ by
% cf Lo67, (loc. cit. Satz 1.5)
$$
\nu([f,v],[g,w])=
[f\sigma(f)^{-1}\sigma(g),w],
$$
and one can show that $G \times_H U$ becomes a reflection space
such that the projection onto $M$ becomes a homomorphism of reflection spaces
(cf.\ [Lo67, Satz 1.5]).
Let us say that then {\it $F$ is a reflection space over the
symmetric space $M$}.
The preceding formula  shows that the reflection $\nu_v$
does not depend on the choice of $v \in U$, i.e., it depends only on the base.
O.\ Loos has shown ([Lo67]) that, conversely, every real finite-dimensional
and connected reflection space can be written in this way as a 
homogeneous bundle over a symmetric space.
In particular, linear representations of $H$ and reflection spaces over
$M$ with linear fibers (``reflection vector bundles over $M$'') correspond to
each other. 

\ssk
Having this in mind, we now consider a symmetric bundle
$\pi:F \to M$ over a homogeneous symmetric space $M=G/H$. As we have just seen,
 $F$ is then also homogeneous, say, $F=L/K$. 
Looking at $H$ as a subgroup of $K$ (see 1.8), we get a linear representation
of $H$ on the fiber $F_o$, and we can  write $F =  G \times_{H} F_o$
 as a homogeneous bundle over the
base $M$. In this way, the functor from symmetric bundles to reflection spaces
corresponds in the homogeneous case to the functor from symmetric 
bundles $F=L/K$ over $M=G/H$ to the associated homogeneous bundle 
$F=G \times_H F_o$. Conversely, we can formulate an
{\it extension problem}:
{\it For which representations $H \to \Gl(V)$ does the homogeneous bundle
$F=G \times_H V$ carry a symmetric bundle structure? If it does, how
many such structures are there?}
%[idea: these are sort of holonomy representations for the connection
%on $F$, i.e.: there is a connection such that the given representation of
%$H$ becomes a holonomy rep.]

\msk
\nin {\bf 1.11. Derivations of symmetric bundles, and vertical automorphisms.}
A {\it derivation} of a symmetric bundle $F$ is a homomorphism of symmetric spaces
$X:M \to F$ which at the same time is a smooth section of $\pi$
 (see [Lo69] for this terminology in case of the tangent bundle).
A {\it vertical automorphism} of a symmetric bundle is an automorphism
$f$ of the symmetric bundle $F$ preserving fibers, i.e., $\pi \circ f =
\pi$. Clearly, $X:= f \circ z:M \to F$ then is a derivation of $F$.
Conversely, if $X$ is a derivation, define
$$
f:F \to F, \quad F_p \ni v \mapsto v + X(p).
$$
Then $f$ is a vertical automorphism: it clearly is smooth, preserves
fibers and is bijective. It is an automorphism: using (SB2),
$$
\eqalign{
\mu(f(v),f(w)) & = \mu(v+X(p),w+X(q)) = \mu ((v,w)+(X(p),X(q))) \cr
&= \mu(v,w) + \mu(X(p),X(q)) =
\mu(v,w) + X \mu(p,q) = f(\mu(v,w)). \cr}
$$
Summing up, vertical automorphisms are the same as derivations.
Moreover, they clearly form a normal subgroup $\Vaut(F)$
in the group $\Aut(F)$,
where composition corresponds to addition of sections.
It follows that the space of derivations is stable under addition;
it is also stable under multiplication by scalars, hence forms a 
vector group.  The same kind of arguments shows that in fact
 we have an exact sequence
$$
0 \to \Vaut(F) \to \Aut(F) \to \Aut(M) \to 1
$$
(which essentially splits if we take transvection groups; 
cf.\ the discussion in 1.9). Now fix a base point $o \in M$; then
the involution given by conjugation with $\sigma_{0_o}$
restricts to $\Vaut(F)$ and thus defines a linear map. Let us write
$$
\Vaut(F) = \Vaut^+(F) \oplus \Vaut^-(F) 
$$
for the corresponding eigenspace decomposition.

\Lemma 1.12. 
$$
\Vaut^+(F) = \{ f \in \Vaut(F) | \, f(0_o)=0_o \},
$$
and the map
$$
\Vaut^-(F) \to F_o, \quad f \mapsto f(0_o) 
$$
is a bijection wit inverse $v \mapsto \sigma_{v \over 2} \sigma_{0_o}$.
 
\Proof. See [Lo69] or
 [Be08], Prop. 5.9, for the proof in the case of the tangent bundle;
the same arguments apply here.
\qed

\sectionheadline
{2. General representations of Lie triple systems}

\Definition 2.1. A {\it Lie triple system (Lts)} 
is a linear space $\m$ over $\K$ together
with a trilinear map $\m \times \m \times \m \to \m$, $(X,Y,Z) \mapsto [X,Y,Z]$ 
such that, writing also $R(X,Y)$ for the endomorphism $[X,Y,\cdot]$,

\item{(LT1)} $R(X,Y)=-R(Y,X)$ (skew-symmetry)
\item{(LT2)} 
$R(X,Y)Z+R(Y,Z)X+R(Z,X)Y=0$ (the Jacobi identity), 
\item{(LT3)} $R(X,Y)$ is a derivation of  the trilinear
product on $\m$, i.e.
$$
R(X,Y)[U,V,W] = [R(X,Y)U,V,W]+[U,R(X,Y)V,W]+[U,V,R(X,Y)W].
$$

\nin For instance, if $(\g,\sigma)$ is a Lie algebra with involution, then
the $-1$-eigenspace $\m$ of $\sigma$ with $[X,Y,Z]:=[[X,Y],Z]$ is a Lts.
Every Lts arises in this way (see Section 3.1 below).
For later use we introduce also the ``middle multiplication operators''
$M(X,Z)Y:=[X,Y,Z]$; then, in presence of (LT1), property (LT2) can be written
in operator form

\ssk
\item{(LT2a)}
$M(X,Z) - M(Z,X)=R(X,Z)$,

\ssk \nin
and similarly, reading (LT3) as an identity of operators,
applied to the variable $W,V$ or $Y$, we get the following
equivalent conditions, respectively

\ssk
\item{(LT3a)}
$[R(X,Y),R(U,V)] =R(R(X,Y)U,V)+R(U,R(X,Y)V)$,
\item{(LT3b)}
$[R(X,Y),M(U,W)] = M(R(X,Y)U,W) + M(U,R(X,Y)W),$
\item{(LT3c)}
$M(X,R(U,V)W) = - M(V,W) \circ M(X,U) + M(U,W) \circ M(X,V) +
R(U,V) \circ M(X,W).$

\msk \nin
{\bf 2.2. The Lie triple system of a symmetric space.}
Let $(M,\mu)$ be a symmetric space with base point $o$.
Consider the tangent bundle $TM$ and write $\g:=\Vaut(TM)$ for the
derivations of the symmetric bundle $TM$,
and $\g=\h \oplus \m$ for the eigenspace decomposition from Lemma 1.12.
One shows that $\g$, seen as a space of vector fields on $M$, is stable
under the Lie bracket and that $\sigma_0$ induces an involution of this
Lie algebra structure. Hence the $-1$-eigenspace $\m$ is a Lts.
Via the bijection $\m \to T_o M$, $X \mapsto X(o)$ from Lemma 1.12 this
Lts structure can be transferred to the tangent space $T_oM$, which, by
definition, is the {\it Lts associated to the pointed symmetric space
$(M,o)$} (cf.\ [Lo69] or [Be08, Chap 5]).
The Lts depends functorially on $M$ and plays a similar role
for symmetric spaces as the Lie algebra for a Lie group. (In particular,
in the real finite dimensional case there is an equivalence of categories
between Lts and connected simply connected spaces with base point, cf.\
[Lo69].)

\msk
\nin {\bf 2.3. The Lts of a symmetric bundle.}
Now assume that $\pi:F \to M$ is a symmetric bundle over
$M$ and  fix a base point $o \in M$ and let 
$\f := T_{0_o} F$. 
Since $F$ is a symmetric space, $\f$ is a Lie
triple system. We wish to describe its structure in more detail.
The differentials of the three involutions $\sigma_0$, $\vartheta_0$ and
$\nu_0$ from Section 1.6
act by automorphisms on the Lts $\f$ (where we write $0$ instead of $0_o$). The
$+1$-eigenspace of $\vartheta$ can be identified with the Lts $\m$, via
the tangent map of the zero-section, and its $-1$-eigenspace is the tangent
space $T_0 (F_o)$ of the fiber $V=F_o$ which we identify with $V$
(in the notation from [Be08] we could also write $\eps V$ for this ``vertical
space'').

\Lemma 2.4. The decomposition
$$
\f = \m \oplus V
$$
has the following properties:

\item{(1)}
$V$ is an ideal of $\f$, i.e., $[x,y,z] \in V$ as soon as one of the 
$x,y,z$ belongs to $V$,
\item{(2)}
$\m$ is a sub-Lts of $\f$,
\item{(3)}
$[V,V,\m]=[V,\m,V]=0$,
\item{(4)}
$V$ is abelian, i.e., $[V,V,V]=0$.

\Proof.
 (1) is clear since $V$ is the kernel of a homomorphism (the differential
$T_0\pi$ 
of $\pi$ at $0$); (2) follows from the fact that $\m$ is the fixed point
space of the horizontal automorphism $T_0(\vartheta_0)=(-1)_F$ (moreover, we see that
$\m$ is isomorphically mapped by $T_0\pi$ onto the Lts of $M$),
and  (3) holds since $\tilde \nu:=T_0(\nu_0)$
is an automorphism and hence
$$
\tilde \nu [v,w,X]=[\tilde \nu(v),\tilde \nu(w),
\tilde \nu(X)]=[v,w,-X]=-[v,w,X],
$$
hence $[v,w,X]$ belongs to $\m$ and thus
to $\m \cap V =0$ since $V$ is an ideal.
Finally, the fiber $F_o$ carries the ``flat'' symmetric space structure of a linear
space and hence
$[V,V,V]=0$. 
\qed

\nin A side-remark: if we adapt the whole set-up to the case of {\it bilinear bundles}
(in the sense of [Be08]) instead of vector bundles, we get essentially the same
properties, the only difference being that (3) does no longer hold:
e.g., for the bilinear bundle $TTM$ over $M$, the tangent 
model is  $\m \oplus (\eps_1 \m \oplus \eps_2 \m \oplus \eps_1 \eps_2 \m)$,
where the term in brackets 
%is the ``$\m$-module'' in that category. 
%This term 
is still an abelian 
Lts, but (3) does no longer hold: in fact, $[\eps_1 \m,\eps_2 \m,\m]$ is non-zero
in general.

\msk
\nin
{\bf 2.5. Representations and modules.}
Let $\m$ be a Lts. An {\it $\m$-module} is a vector space $V$ such that
the direct sum $\f:=\m \oplus V$ carries the structure of an Lts satisfying
the properties from the preceding lemma. More explicitly, this means,
by decomposing
$$
[X\oplus u,Y \oplus v,Z\oplus w] =
[X,Y,Z] \oplus (r(X,Y)w + m(X,Z)v - m(Y,Z)u + [u,v,w]) 
\eqno (2.1)
$$
that we are given two trilinear maps $r$ and $m$ 
$$
\eqalign{
r: \m \times \m \to \End(V), & \quad r(X,Y)=[X,Y,\cdot], \cr
m: \m  \times \m \to \End(V), & \quad m(X,Z)=[X,\cdot,Y]=
- [\cdot,X,Y] \cr}
\eqno (2.2)
$$
satisfying the properties given by the following lemma:

\Lemma 2.6.
For any Lts $(\m,R)$, the
space $\m \oplus V$ with a triple bracket given by
{\rm (2.1)} is a Lts if and only if $r$ and $m$ satisfy the following
relations:

\item{(R1)}
$r(X,Y) = - r(Y,X)$,
\item{(R2)}
$m(X,Z) - m(Z,X)=r(X,Z)$
\item{(R3)} $r(R(X,Y) U \otimes V+U \otimes R(X,Y)V)=[r(X,Y), r(U,V)]$,
\item{ }
$m(R(X,Y)  U \otimes V+ U \otimes R(X,Y)V)=[r(X,Y),m(U,V)]$,
\item{(R4)}
$m(X,R(U,V)W)-r(U,V) \circ m(X,W) =
 m(U,W) \circ m(X,V)   - m(V,W) \circ m(X,U)$.
 
\Proof.
We have to show that (LT1) -- (LT3) for $\m \oplus V$ are
equivalent to (R1) -- (R4):
first of all, we note that a bracket is zero if more than
one of the three arguments belongs to $V$.
Now, (LT1) is equivalent to (R1) if both arguments belong to
$\m$ and holds by (2.10) if one is in $V$ and the other in $\m$.
Next,
(LT2) is an identity in three variables. We may assume that
two variables, say $X$ and $Z$, belong to $\m$, and write (LT2) in its 
operator form (LT2a). Thus we see that (LT2) is equivalent to (R2).
Finally, (LT3) is an identity in 5
variables. In order to get a non-trivial identity, we can assume
that at least four of them belong to $\m$.
We then write (LT3) in operator form (identities (LT3a,b,c) from
Section 2.1), and see that
(LT3) is equivalent to (R3) and (R4), thus proving our claim.
\qed

\nin Note that,
in view of (R3),  identity (R4) is equivalent to
the following identity:

\ssk
\item{(R4')}
$m(R(U,V)X,W)-m(X,W) \circ r(U,V) =
 m(V,W) \circ m(X,U) -m(U,W) \circ m(X,V).$

\ssk
\nin
Condition (R4) can be rephrased by saying that the
operator $R(X,v)$ defined by $R(X,v)Y=m(X,Y)v$  belongs 
to the space of {\it derivations from $\m$ into $V$},
$$
\eqalign{
\Der(\m,V) &  =\{ D:\m \to V | \, \forall X,Y,Z \in \m : \cr
& \quad \quad
D R(X,Y,Z)=r(X,Y)DZ + m(X,Z)DY - m(Y,DX)Z \}.
\cr}
$$

\Definition 2.7.
A {\it general representation of a Lie triple system}
$\m$ in a unital associative algebra $A$
is given by two bilinear maps
$$
\eqalign{
r:\m \times \m \to A, &
 \quad (X,Y) \mapsto r(X,Y), \cr
m:\m  \times \m \to A, &
\quad (X,Z) \mapsto
m(X,Z) \cr}
$$
% (equivalently,
% $$r:\q \times \q \to A, \quad m:\q \times \q \to A,$$
%with $A$ the associative algebra $\End(V)$)
such that (R1) - (R4) hold (where $\circ$ has to be interpreted as the product in $A$
and the bracket is the Lie bracket in $A$).
If $A = \End(V)$ is the endomorphism algebra of a vector space, 
we  say  that $V$ is an {\it $\m$-module}.
{\it Homomorphisms of $\m$-modules} are defined
in the obvious way, thus turning $\m$-modules into a
category. Given an $\m$-module $V$, 
the Lts $\tilde \m = \m \oplus V$ 
whith bracket defined by (2.1) is called the {\it split null extension
of $\m$ by the module $V$}.
It is fairly obvious that the split null extension depends
functorially on the $\m$-module $V$.

\Example 2.8.  
(Regular representation.)
For any Lts $\m$, consider its ``extension by dual numbers'', i.e.,
let $\K[\eps]=\K[X]/(X^2)=\K \oplus \eps \K$, $\eps^2=0$ (ring of dual numbers over
$\K$), and 
$$
\tilde \m = \m \otimes_\K \K[\eps]= \m \oplus \epsilon \m ,
$$
with the $\epsilon$-trilinear extension of the bracket from $\m$:
$$
[X+\epsilon X',Y+\epsilon Y',Z + \epsilon Z']=
[X,Y,Z] + \epsilon ([X,Y,Z']+[X,Y',Z]+[X',Y,Z]).
\eqno (2.3)
$$
This is nothing but the split null extension of $\m$ by the {\it regular
representation}, which by definition is given by $V=\m$ and
$$
r(X,Y)=R(X,Y):\m \to \m, \quad \quad
m(X,Y)=M(X,Y):\m \to \m .
$$
If $M$ is a symmetric space, then  
the Lts of the tangent bundle $TM$ is precisely $\m \oplus \eps \m$ (cf.\ [Be08]).
Hence the regular representation corresponds to the tangent bundle of $M$.

\ssk \nin {\bf 2.9. The extension problem revisited.}
The forgetful functor associating to a symmetric
vector  bundle its underlying reflection space corresponds
to the forgetful functor $(r,m) \mapsto r$. 
%associating to
%a $\q$-module its underlying $\h$-module.
Namely, if $\h$ is the image of the skew-symmetric map
$$
\m \otimes \m \to \End(\m), \quad X \otimes Y \mapsto R(X,Y),
$$
then (LT3) implies that $\h$ is a Lie algebra, and for any representation
$\rho:\h \to \gl(V)$, we may define $r(X,Y):=\rho(R(X,Y))$;  then the first
relation of (R3) is equivalent to $\rho$ being a representation.
Thus we get the infinitesimal version of a homogeneous vector bundle.
Now the problem of finding a compatible symmetric vector bundle structure corresponds
to finding the second component $m$ such that $(r,m)$ defines a representation of
$\m$.

%
%\nin
%At present we have no general answer to these questions.
%However, it seems clear that the answer should be formulated
%in terms of some cohomology theory. Here the most interesting
%references are [Ha61] and [McK71].

\sectionheadline{3. Reconstruction}

We have shown that a representation of an Lts is the derived version
of a symmetric  bundle. Conversely, can one reconstruct a symmetric
bundle from a representation of an Lts? 
As a first step, it is always possible to recover Lie algebras from
Lie triple systems, and certain Lie algebra representations from
Lie triple representations. 
The second step is then to lift these constructions to the space level:
here we have to make assumptions on the base field and on the topological nature of $M$.

\subheadline{From Lie triple systems to Lie algebras with involution}

\nin
{\bf 3.1. Lie triple systems and $\Z/2\Z$-graded Lie algebras: the standard
imbedding.}
Every Lie algebra $\g$ together with an involution $\sigma$ gives rise to
a Lts $\m = \g^{-\sigma}$, equipped with the
triple Lie bracket $[[X,Y],Z]$.
Conversely, every Lts $\m$ can be obtained in this way:
let $\h$ be the subalgebra of the algebra of derivations of $\m$ generated
by the endomorphisms $R(x,y)$, $x,y \in \m$. Then the space
$$
\g:= \h \oplus \m
$$
carries a Lie bracket given by
$[(D,X),(D',X')]:=([D,D'] + R(X,X'),DX'-D'X)$.
This Lie algebra, called the {\it standard imbedding of the Lts $\m$},
does in general not depend functorially on $\m$ -- see [HP02] and
[Sm05] for a detailed study of functorial properties related to this
and other constructions.
Note that, in terms of the Lie algebra $\g$, we can write
$$
R(X,Y)=\ad[X,Y]|_\m, \quad \quad
M(X,Z)=\ad(Z) \circ \ad(X) |_\m.
\eqno (3.1)
$$

\msk
\nin {\bf 3.2. $\g$-modules with involution.}
Assume $(\g,\sigma)$ is a Lie algebra with involution.
A representation $\rho:\g \to \gl(W)$ is called a
 {\it $(\g,\sigma)$-module with involution} if $W$ is equipped with
a direct sum decomposition $W=W^+ \oplus W^-$ which is compatible with
$\sigma$ in the sense that
$$
\matrix{
\g & {\buildrel \rho \over \longrightarrow} & \gl(W) \cr
\sigma\downarrow \phantom{\sigma}& &\phantom{\tau_*} \downarrow \tau_*
 \cr
\g & {\buildrel \rho \over \longrightarrow} & \gl(W) \cr}
$$
commutes, where $\tau \in \Gl(W)$ is the identity on $W^+$ and $-1$ on $W^-$, and
$\tau_*(X)=\tau X \tau$.

\Lemma 3.3.
Let $\m$ be a Lts and $\g$ its standard imbedding.
There exists a bijection between $(\g,\sigma)$-modules with involution
and $\m$-modules.

\Proof.
Given a $(\g,\sigma)$-module with involution $(W,\tau)$, we first form
the semidirect product $\b:=\g \ltimes W$. This is a Lie algebra carrying
an involution given by $\sigma \times \tau$. Its $-1$-eigenspace
$\m \oplus W^-$ is a Lts satisfying the relations from Lemma 2.4, and
hence is the split null extension corresponding to an $\m$-module $W^-$.

Given an $\m$-module $V$, we construct first the split null extension
$\m \oplus V$ and then its standard imbedding $\b = (\m \oplus V) \oplus
[\m \oplus V,\m \oplus V]$. Then 
$$
W:= W^+ \oplus W^- := [V,\m] \oplus V
$$
is a $\g$-module with involution. 
%(Elements of $[V,\m]$ act as nilpotent
%operators $\m \oplus V \to V \to 0$.) 
\qed

\nin Again, the correspondence set up by the lemma is functorial in one direction
but not in the other -- see [HP02] for this issue.

\subheadline{From modules to bundles}

\Theorem 3.4. 
Let $\K=\R$ and $M$ be a finite-dimensional connected simply connected symmetric space
with base point $o$. Let $\m$ be its
 associated Lts and $\g$ its standard imbedding, with involution $\sigma$.
 Then the following objects are in
one-to-one correspondence:

\item{(1)}
(finite dimensional) symmetric vector bundles over $M$,
\item{(2)}
(finite dimensional) $(\g,\sigma)$-modules with involution,
\item{(3)}
(finite dimensional) $\m$-modules.
\vskip 0mm
\ssk \nin
The bijection between {\rm (1)} and {\rm (3)} is an equivalence
of categories.

\Proof. 
We have already seen how to go from (1) to (3), and that (2) and (3)
are in bijection. Let us give a construction from (2) and (3) to (1):
let $\f=\m \oplus V$ be the split null extension coming from an
$\m$-module $V$ and $\b=\f \oplus [\f,\f]$ the standard imbedding of $\f$
and let $W$ the corresponding $\g$-module with involution.
Let $G$ be the simply connected covering of the transvection group $G(M)$ and
write $M=G/H$.
 Then the representation of $\g$ on $W$ integrates to a
representation of $G$ on $W$. Let
$$
B:= G \ltimes W, \quad \quad K:= H \ltimes W^+, \quad \quad F:= B/K.
$$
We claim that

\item{(i)} $F$ is a vector bundle over $M$, isomorphic to the
homogeneous bundle $G\times_H W^-$, and 
\item{(ii)} $F$ carries the structure of
a symmetric bundle over $M$.

\ssk \nin
Proof of (i):
first of all,
$$
GW/HW^+ \to G \times_H (W/W^+), \quad gw/HW^+ \mapsto [(g,w / W^+)]
$$
is a well-defined bijection.
Since $W/W^+ = W^-$, this proves the first claim.

Proof of (ii):
the Lie algebra  $\k$ of $K$ is the fixed point space of an involution of $\b$,
and hence $F=B/K$ is a symmetric space. Its Lts is $\f$. The projection map
$F \to M$ has as differential the projection from $\f$ to $\m$ and hence is
a homomorphism of symmetric spaces.  

Let us show that the structure map $F_p \oplus F_q \to F_{\mu(p,q)}$ is linear.
Since we already know that $F$ is a homogeneous $G$-bundle, we may assume
that $p=o$ is the base point. 
Now we proceed in two steps:

\ssk \nin
(a) we show that $\sigma_{0_o}:F_{\sigma(q)} \to F_q$ is linear.
In fact, here we use that $W$ is a $G$-module with involution, i.e.,
$\sigma(g)w=\tau g \tau(w)$:
$$
\eqalign{
\sigma([(g,w)]) & = \sigma(gw) /HW^+ = \sigma(g)\sigma(w)/HW^+ \cr
& = \sigma(g)(-w)/HW^+ = \tau(g(w))/HW^+ = - g(w)/HW^+ =
[(g,-w))].\cr}
$$
Thus this map is described by $w \mapsto -w$ and thus is linear.
%
%[Other argument: since the involution on the fiber $F_o = W^-$ is just minus
%the identity, which is linear, we may replace without loss $H$ and $G$ by the
%groups generated by $-\id$ and $H$, resp. $G$, then these groups act linearly
%on fibers.]

\ssk \nin
(b)
Since $\sigma_{0_o}:F_{\sigma(q)} \to F_q$ is a linear bijection,
$F_o \oplus F_q \to F_{\sigma(q)}$ is linear (and well-defined) iff so is 
the map $F_o \oplus F_q \to F_q$, $(u,w) \mapsto \sigma_{0_o} \sigma_{u}(w)$.
But the last map is the same as $(u,w) \mapsto (-2u).w$ (the point stands for
the action of $W^-$ on $F$; recall that  in every symmetric
space $\sigma_o \sigma_{g.o}=\sigma(g)g^{-1}$). 
Summing up,
it suffices now to show that the map
$$
W^- \times F_q \to F_q, \quad (u,z) \mapsto u.z
$$
is well-defined and linear. Proof of this:
let $u \in W^-$ and $q=gH$, $[(g,w)] \in F_q$ with $w \in W^-=W/W^+$. Then
$$
\eqalign{
u.[(g,w)] & = u.gw /HW^+ = gg^{-1}ugw/HW^+ \cr
& = g(\rho(g)^{-1}u)w/HW^+ = [(g,w+\pr_-(\rho(g)^{-1}u))].\cr}
$$
Thus our map is described by $(u,w) \mapsto w+\alpha(u)$
with a linear map $\alpha \in \End(W^-)$ that depends on $g$, and hence is linear,
proving claim (ii).
%
%(reconstruction proof in the setting of [Did06]
% was simpler because one just had to show that the horizontal map
%is an automorphism!)

\ssk
Finally, the fact that homomorphisms in the categories defined by (1) and (3)
correspond to each other follows from the corresponding fact for 
(connected simply connected) symmetric spaces and Lie triple systems, see [Lo69].
\qed

\sectionheadline
{4.  Linear algebra and representations}

So far we do not know any representations other than the regular one
and the trivial ones. In the following we discuss the standard linear
algebra constructions producing new representations from old ones:

\msk
\nin {\bf 4.1. Direct sums.}
Clearly, if 
 $(V,r_1,m_1)$ and $(W,r_2,m_2)$ are $\m$-modules, then
$(V \oplus W, r_1 \oplus r_2, m_1 \oplus m_2)$ is again
a general representation.
Correspondingly,
 the direct sum of symmetric bundles can be turned into
a symmetric bundle. 

\msk
\nin {\bf 4.2. The dual representation.}
If $(V,r,m)$ is a $\m$-module, then
the dual space $V^*$ can be turned into an $\m$-module by putting
$$
r^*(X,Y):=-r(X,Y)^*=r(Y,X)^*, \quad
m^*(X,Y):=m(Y,X)^*,
\eqno (4.1)
$$
where $A^*:V^* \to V^*$, $\phi \mapsto \phi\circ A$ is the
dual operator of an operator $A \in \End(V)$.
In fact,
the properties (R1) - (R3) for $(V^*,r^*,m^*)$
are easily verified; for (R4)
note that (R4') written out for the dual is precisely (R4).

Equivalently:
if $(\rho,\g,V,\tau)$ is a $(\g,\sigma)$-module with involution, one
 verifies  that the dual module,
$\rho^*(X)=-\rho(X)^*$,  also is a module with involution $\tau^*$.
It follows that
$$
m^*(X,Y)=\rho^*(X)\rho^*(Y)=
(-\rho(X)^*)(-\rho(Y)^*)=(\rho(Y)\rho(X))^*=m(Y,X)^*
$$
leading to Formula (4.1).

In particular, in the finite dimensional real case, invoking Theorem 3.4,
the dual of the regular
representation corresponds to the cotangent bundle $T^* M$ which thus again
carries a symmetric bundle structure. (If $M=G/H$, then we may also write
$TM=T^*G/T^*H$ where $T^*G=G \ltimes \g^*$.)
It remains intriguing  that there seems to be no really intrinsic
construction of this symmetric space structure on $T^*M$.
For this reason we cannot affirm that (in cases where
 a reasonable topological dual $\m^*$ of $\m$ exists),
 in the infinite dimensional case
or over other base fields than $\R$ or $\C$, $T^*M$ is again a symmetric space. 

%Note: on the level of bundles we cannot make a similar construction!!
%However, in the finite-dimensional real case, as topological spaces,
%bundles and dual bundles are always isomorphic (there is always a
%Riemannian metric $F \to F^*$). Idea: just declare that $F^*=F$ as set
%and has the same multiplication map $F \times F \to F$; but declare that
%in $Diff(F)$ multiplication is reversed. This does not affect the symmetric
%space axioms, and we still get the same symmetric space.
%However, the identity map is no longer equivariant with respect to the
%action of the respective transvection groups. In fact, we will et
%minus-signs in the appropriate places.

\msk
\nin {\bf 4.3. A duality principle.}
Note that, for finite dimensional modules over a field,
$V$ is the dual of its dual module $V^*$.
More generally, we can define for any general representation
of $\q$ in an algebra $A$ its {\it opposite representation} in
the algebra $A^{opp}$ by putting
$$
r^{opp}(X,Y)=r(Y,X), \quad \quad m^{opp}(X,Y)=m(Y,X);
$$
as above it is seen that this is again a representation.
As an application of these remarks  we get a
{\it duality principle} similar as the one for
Jordan pairs formulated by O. Loos (cf. [Lo75]):

\Proposition 4.4.
If $I$ is an identity in $R(X,Y)$ and $M(U,V)$ valid
for all  Lie triples over $\R$, then its dual identity $I^*$,
obtained by replacing $R(X,Y)$ by $R(Y,X)$
and $M(U,V)$ by $M(V,U)$ and reversing the order of
all factors, is also valid for all
 Lie triples over $\R$.

\Proof.
If $I$ is valid for all Lts, then it is also valid
for all split null extensions obtained from  representations
and hence the corresponding identity, with $R(X,Y)$ replaced
by $r(X,Y)$ and $M(X,Y)$ by $m(X,Y)$, is valid for all 
representations .
 Since the set of all  representations
is the same as the set of all opposite representations, and since
the opposite functor changes order of factors and 
order of arguments,
we see that $I^*$ is valid for all representations.
In particular, it is valid for the regular representation
and hence holds in $\q$.
\qed

\nin
For instance, identities (R4) and (R4') (cf.\ Lemma 2.6)
 are dual in the sense of the proposition.
We don't know about any application of Proposition 4.4;
however, one may note that the original definition of
Lie triple systems by N. Jacobson in [Jac51] as well
as the exposition by Lister [Li52] are based on a set
of five identities, among which two identities  are
equivalent to each other by the duality principle
-- they correspond to (LT3c) and its dual
identity.

\msk \nin
{\bf 4.5. Tensor products.}
The tensor product $F_1 \otimes F_2$ of two symmetric vector bundles is
in general no longer a symmetric vector bundle:
let $A:=(F_1)_o$, $B=(F_2)_o$ be the two fibers in question, regarded as
$\m$-modules.
Extend $A$ and $B$ to $(\g,\sigma)$-modules with involution,
$V=V_+ \oplus V_-$, $A = V_-$, $V_+=\Der(\m,A)$,
$W=W_+ \oplus W_-$, $B=W_-$, $W_+ = \Der(\m,B)$.
It is easily verified that then
$(V \otimes W,\tau_V \otimes \tau_W)$ is again a
$(\g,\sigma)$-module with involution.
Now, the minus-part in
$$
V \otimes W = (V_+ \otimes W_+ \oplus V_- \otimes   W_-) \oplus
(V_+ \otimes W_- \oplus V_- \otimes   W_+),
$$
 is
$$
A \odot  B := A \otimes \Der(\m,B) \oplus B \otimes \Der(\m,A)
$$
which therefore is another $\m$-module, replacing the
ordinary tensor product $A \otimes B$.
The corresponding definition of the maps $r$ is
$$
r(X,Y)(a \otimes D )=
r_A(X,Y)a \otimes D + a \otimes (r_B(X,Y) \circ D - D \circ R(X,Y),
$$
and there is a similar expression for the $m$-components and
 for $r(X,Y)(b \otimes D')$, $m(X,Y)(b\otimes D')$.
It is obvious that the operation $\odot$ is compatible with
direct sums, and it also associative (in the same sense as the usual
tensor product): the minus-part both in $U \otimes (V \otimes W)$
and in $(U \otimes V) \otimes W$ is 
$$
U_- \odot V_- \odot W_- = \bigoplus_{ijk=-1} U_i \otimes V_i \otimes W_i
$$
where $i,j,k \in \{ \pm 1 \}$ and
with $U_+ =\Der(\q,U_-)$, etc.

\msk
\nin
{\bf 4.6. Hom-bundles.}
If $A$ and $B$ are $\m$-modules, then $\Hom(A,B)$
is in general not an $\m$-module, but we can use the same
construction as in 4.5 to see that
$$
\Hom(A,\Der(\m,B)) \oplus \Hom(B,\Der(\m,A))
$$
is again a $\m$-module.

\msk
\nin {\bf 4.7. Universal bundles.}
Various definitions of {\it universal} or {\it enveloping algebras}, resp.\
{\it universal representations}, attached to triple systems or algebras with involution,
can be given -- see [Ha61], [Lo73], [Lo75], [MoPe06].
These objects should correspond to certain ``universal bundles'' over a given symmetric
space $M$. We intend to investigate such questions elsewhere.

\sectionheadline{5. The canonical connection of a symmetric bundle}

We now study in more detail the differential geometric aspects of
symmetric bundles. It is immediately clear from Section 1.6 that a
symmetric bundle $F$ carries a fiber bundle connection in the general
sense of Ehresmann, i.e., there is a distribution of {\it horizontal
subspaces}, complementary to the vertical subspaces
 $V_u=T_u(F_p)$ (tangent spaces of the
fiber) -- namely, as horizontal subspace take the fixed point spaces of
the differentials of the horizontal symmetry $\vartheta_u$,
$$
H_u= \{ v \in T_u F | \, \vartheta_u (v)=v \} = Q({u \over 2})(T_p M) =
\sigma_{u\over 2}(T_pM).
$$
In the sequel, we show that this Ehresmann connection is indeed a {\it
linear connection} (in the general sense defined in [Be08], which in the
real case amounts to the usual definitions), and that in general it has
non-vanishing curvature, so that the distribution $(H_u)_{u \in F}$ 
is in general not integrable.
 
\Theorem 5.1.
Let $\pi:F \to M$ be a symmetric bundle over $M$. Then there exists a unique
linear connection on the vector bundle $F$ which is invariant under all
symmetries $\sigma_x$, $x \in M$. 
%A symmetric vector bundle carries a canonical Ehresmann connection.

\Proof.
The proof is similar the one of [Be08, Theorem 26.3] and
will therefore not be spelled out here in full detail.
The uniqueness statement is proved by observing that the difference $A=L_1 - L_2$
of two linear connections on $F$ is a tensor field such that
$A_p:T_p M \times F_p \to F_p$ is bilinear; if both $L_1$ and $L_2$ are invariant
under $\sigma_p$, then $A_p(-v,-w)=-A_p(w)$ and hence $A_p=0$.
It follows that $L_1=L_2$.
The main argument for the proof of existence
consists in proving that the fibers of the bundle
$TF$ over $M$ are abelian symmetric spaces (for this one has to analyze
the map $T\mu:TF \times TF \to TF$ in the fiber over a point $p \in M$
in the same way as the corresponding map $(TTM \times TTM)_p \to (TTM)_p$ was analyzed
in [Be08, Lemma 26.4]); then one concludes by general arguments that
the fibers of $TF$ over $M$ carry canonically the structure of a linear space
which is bilinearly related to all linear structures induced by bundle charts.
By definition, this is what we call a linear connection on $F$. 
\qed

%
%Achtung: ein ALLGEMEINER Ehresmann-Zhg muss nicht unbedingt bilinearly related mit
%den Kartenzhgen sein [ie., es muss kein linearer Zhg sein. Anders gesagt:
%der Buendelatlas muss nicht derselbe sein -- Buendel koennen abstrakt
%isomorph sein aber verschiedene Buendelatlanten haben!z.b. Geradenbuendel
%mit $x \mapsto x^3$ in jeder Faser...]. 
% Unsere Behauptung schliesst ein, dass dieser Ehresmann-Zhg bilin.related ist,
% und dafuer ist loc cit Lemma 26.4 unumgaegnlich, wie mir scheint!
%Hier schonmal eine Konsequenz: falls $F=TM$, so liefert jedes symmetrische
%denselben Zhg, also charakterisiert dieser NICHT eindeutig das symmetrische
%Buendel!
%
%24.6: oder vielleicht ganz einfach so: es stimmt doch alles wie geschrieben,
%und die beiden Zhge auf $TM$ sind ganz einfach dieselben, dh. der Zhg eines
%symmetrischen Buendels charakterisiert dieses Buendel eben nicht.
%Dies tut erst der affine Zhg auf $F$, d.h. der Zhg auf $TF$.
%Und in der Tat, auf $TM$ gibt es ein $G$-invariantes Algebrenfeld: das
%JT-Produkt $T(\cdot,y,\cdot)$, wobei $y$ das Element in der Faser beschreibt
%(und auf dem Nullschnitt verschwindet das Feld, wie erwartet).

\ssk \nin {\bf 5.2. Extension problem: on uniqueness.}
When $F=TM$ is the tangent bundle, with its canonical symmetric bundle structure,
the connection defined by the preceding theorem is precisely the {\it canonical
connection of the symmetric space $M$}, cf.\ [Be08], [Lo69]. 
We will see in the next chapter that the abstract bundle $TM$ may carry several
different (non equivalent) symmetric bundle structures over $M$. The uniqueness
statement of the theorem shows that they all lead to the same linear connection
on $TM$ over $M$. Therefore the symmetric bundle structure is {\it not} uniquely
determined by the linear connection from Theorem 5.1.
Only at second order, by considering connections on $TTM$ over
$TM$, one is able to distinguish two symmetric bundle structures on $TM$.
%
% investigate this further: see remarks below!

\Theorem 5.3. Let  $\Omega$ be the curvature tensor of
 the linear connection on the symmetric bundle $F$ over $M$ defined in the preceding theorem.
Then $\Omega$
is given by the $r$-component of the corresponding Lts-representation:
for all sections $X,Y$ of $TM$ and sections $\zeta$ of $F$,
$$
\Omega(X,Y)\zeta = r(X,Y)\zeta
$$
(i.e., with respect to an arbitrary base point $o \in M$, for all $v,w \in T_oM$ and $z \in F_o$,
$\Omega_o(v,w)z=r(v,w)z$).

\Proof.
In case $F=TM$, where $\Omega$ is the curvature of the
 canonical connection of $M$, it is well known that
$\Omega$ is (possibly up to a sign, which is a matter of convention)
 given by the Lie triple system of $M$, i.e.\
$\Omega(X,Y)Z = R(X,Y)Z= [X,Y,Z]$,
see [Lo69], [KoNo69], [Be00] and [Be08] for three  different proofs.
It is inevitable to go into third-order calculations, and therefore
none of these proofs is really short. 
Theorem 5.3 generalizes this result, and it can also be proved in several different ways.
We will here just present the
basic ideas and refer the reader to the above references for details.

(a) Approaches using sections.  
Recall that the curvature $\Omega$ may be defined by
$$
\Omega(X,Y)=[X,Y]_h - [X_h,Y_h],
$$
where $X_h:F \to TF$ is the horizontal lift of a vector field $X:M \to TM$. 
According
to the definition of the horizontal space $H_u$ given above,
the horizontal lift of a vector field $X$ is given
by, for $u \in F_x$,
$$
X_h(u) = \sigma_{u\over 2} \sigma_{0_x}(X(x)).
$$
By homogenity, it suffices to calculate $\Omega_o$, the value at the base
point, and then it is enough to plug in the vector fields
$X=\tilde v$, $Y = \tilde w \in \m$ having value $X(o)=v$,
$Y=w$ for $v,w \in T_oM$. This already implies that $[X,Y]_o=0$, and so we
are left with calculating $[\tilde v_h,\tilde w_h]_u$. 
This can be done by a calculation in a chart, corresponding essentially to [Be08],
Lemma 26.4;
the outcome is $[v,w,u]$, as expected.

(b) Approaches using higher order tangent bundles.
We analyze the structure of the bundles $TF$ and $TTF$ in exactly the same
way as done in [Be08, Chapter 27] for the case of the tangent bundle: as
mentioned in the proof of the preceding theorem, $TF$ has abelian fibers
defining the canonical connection, and $TTF$ is non-abelian, leading to an
intrinsic description of the curvature in terms of the symmetric space structure
of the fibers. The conclusion is the same as with $T^3 M$:
$\Omega(u,v)w=[u,v,w]$ for $u,v,w$ belonging to the three ``axes''
$\eps_1 T_pM$, $\eps_2 T_pM$, $F_p$ of $TTF$.
\qed

\ssk \nin {\bf 5.4. Extension problem: on existence.}
Via Theorem 5.3, we can relate the extension problem  to a problem on
{\it holonomy representations}:
assume $M=G/H$ is a homogeneous symmetric space and assume
given a representation $\rho:H \to \Gl(V)$; then
the associated bundle $F=G \times_H V$ carries a tensor field of curvature type
 $r:\m \wedge \m \to \h \to \gl(V)$, coming from the derived representation
$\dot \rho:\h \to \gl(V)$.
 Can we find a  connection (coming from a symmetric bundle structure on $F$)
such that $r$ is its curvature, i.e., such that
this representation becomes a holonomy representation?
In case of the tangent bundle, $F=TM$, the answer clearly is positive,
since $H$ {\it is} the holonomy group of the canonical connection
on the tangent bundle.

\sectionheadline{6. Symmetric structures on the tangent bundle}

In this chapter we will show that, for a rather big class of symmetric
spaces $M$, the tangent bundle $TM$ carries (at least) two different
symmetric structures with the same underlying reflection space structure.
For instance, this is the case for the general linear group, seen as
symmetric space.

\msk
\nin {\bf 6.1. Example: the general linear group.} 
The tangent bundle of the group $\Gl(n,\K)$ can be identified with the 
general linear group over the dual numbers $\K[\eps]$,
$$
T\Gl(n,\K) =  \Gl(n,\K[\epsilon]) = 
\{ g + \eps X | g \in \Gl(n,\K), X \in M(n,\K) \}
$$
with $\epsilon$-bilinear multiplication
$
(g + \epsilon X)(h + \epsilon Y)=
gh + \epsilon(Xh+gY)
$.
The canonical symmetric space structures on $\Gl(n,\K)$ and on its tangent
bundle are given by the product map $\mu(g,h)=gh^{-1}g$ (see Example 1.3 (1)).

\Theorem 6.2. The vector bundle $T \Gl(n,\K)$ admits a second symmetric bundle
structure isomorphic to the  homogeneous symmetric space
$$
L/K= \Gl(n,\D)/\Gl(n,\K[\eps]),
$$
where $\D$ is the non-commutative ring of ``degenerate quaternions over $\K$'', 
$$
\D = \Big\{ \pmatrix{ a & b \cr \overline b & \overline a \cr} | \,
a,b \in \K[\eps] \Big\}
$$
with group involution of $L$ induced by conjugation of $\D$ with respect
to its subalgebra of diagonal matrices.

\Proof.
Before coming to  details of the calculation, let us give a heuristic argument:
the ``Hermitian complexification'' of the group $\Gl(n,\R)$ is the symmetric space
$M_{h\C}=\Gl(2n,\R)/\Gl(n,\C)$ of complex structures on $\R^{2n}$ (see [Be00, Ch.\ IV]).
In principle, in the present context we should have to replace 
complex structures ($I^2 = - \id$) by ``infinitesimal structures'' ($E^2 =0$),
but this attempt fails since $E$ is not invertible. However, 
changing  the point of view by considering a complex structure on $\R^{2n}$
rather as a ``$\C^n$-form of the algebra $(M(2,2;\R))^n$'', and viewing $\Gl(2n,\R)$
rather as $\Gl(n,M(2,2;\R))$, the suitably reformulated arguments carry over from the 
Cayley-Dickson extension $\C \subset M(2,2;\R)$ to the ``degenerate Cayley-Dickson
extension''  $\R[\eps] \subset \D$ (and in fact to any extension \`a la Cayley-Dickson
of a commutative ring with non-trivial involution, see below). -- 
In the following calculations we need the  matrices
$$
F_n= \pmatrix{0 & \1_n  \cr  \1_n & 0 \cr}, \quad
I_{n,n} = \pmatrix{\1_n & 0 \cr 0 & -\1_n \cr}, \quad
R_n:= \pmatrix{\1_n & \1_n \cr -\1_n & \1_n \cr}.
$$
Let us abbreviate $\A := \K[\eps]$; then 
$$
\D = \{ X \in M(2,2;\A) | \, F_1 X F_1 = \overline X \}
$$
where $\overline{A + \eps B}=A - \eps B$. We 
call the map
 $\tau:\D \to \D$, $X \mapsto I_{1,1}XI_{1,1}$
an {\it $\A$-form of $\D$}; this is justified by the fact that
the fixed ring $\D^\tau$ is the ring of diagonal matrices in $\D$, which
is isomorphic to $\A$, and that $\tau$ anticommutes with the
``structure map'' $f:\D \to \D$, $X \mapsto FX$.
For the corresponding $\K$-linear
maps $\D^n \to \D^n$ we will again write $\tau$ and $f$ instead
of $\tau^n$ and $f^n$. 
%
%[The group $\Gl(n,\D)$ acts transitively on the space of 
%$\A$-forms of $\D$, via $(g,c) \mapsto gcg^{-1}$.]
%
The group
$$
\eqalign{
\Gl(n,\D)& = \{ g \in \Gl(2n,\A) | \, F_n g F_n = \overline g \} \cr & =
\Big\{ g=\pmatrix{a & b \cr \overline b & \overline a \cr} | \,
a,b \in M(n,n;\A), \det(g) \in \K^\times \Big\} \cr}
$$
 acts on the space $\End_\K(\D^n)$ by conjugation.
The stabilizer of $\tau$ is given by $b=0$, i.e., it is the group $\Gl(n,\A)$.
Thus  the $\Gl(n,\D)$-orbit of $\tau$ is a homogeneous symmetric space:
$$
{\cal O} := \Gl(n,\D).\tau \cong \Gl(n,\D)/\Gl(n,\A).
$$
We will show now that $\cal O$ is a vector bundle over $\Gl(n,\K)$, and that
this vector bundle is isomorphic to the homogeneous bundle $\Gl(n,\A)$ over $\Gl(n,\K)$.
To this end, observe that the group $\Gl(n,\D)$
contains a subgroup isomorphic to $\Gl(n,\K) \times \Gl(n,\K)$, namely
the group of matrices of the form
$$
\pmatrix{g+h & h-g \cr h-g & g+h \cr} = R_n \pmatrix{g & 0 \cr 0 & h \cr} R_n^{-1}
$$
with $g,h \in \Gl(n,\K)$.
Now let this subgroup act on $\tau$ (whose matrix is $I_{n,n}$):
$$
R_n \pmatrix{g & 0 \cr 0 & h \cr} R_n^{-1} I_{n,n} 
R_n \pmatrix{g^{-1} & 0 \cr 0 & h^{-1} \cr} R_n^{-1} =
R_n \pmatrix{0 & gh^{-1} \cr hg^{-1} & 0 \cr} R_n^{-1}.
$$
The stabilizer of $\tau$ is gotten by taking $g=h$, and so
the orbit of $\tau$ under this group is isomorphic to the
symmetric space $\Gl(n,\K) \times \Gl(n,\K)/diag \cong \Gl(n,\K)$
(group case). 
On the other hand, $\Gl(n,\D)$ contains the abelian
normal subgroup of matrices of the form
$$
\1_{2n} + \eps \pmatrix{X & Y \cr -Y & -X \cr}
$$
with $X,Y \in M(n,n;\A)$, and  $\Gl(n\D)$ is a semidirect product of 
these two subgroups. It follows that $\cal O$ is a vector bundle over 
the orbit $\Gl(n,\K) \times \Gl(n,\K)/diag \cong \Gl(n,\K)$. Let us determine
the fiber over the base point $F_n$.
 Since
$$
(\1 + \eps Z )F_n (\1 - \eps Z)=F_n + \eps (ZF_n - F_nZ),
$$
the stabilizer is gotten by $X=0$, whereas for $Y=0$ we get the fiber of
$\cal O$ over the base point:
$$
\pmatrix{\1  & \eps X \cr - \eps X & \1 \cr}, \quad X \in M(n,n;\K).
$$
Thus the fiber is isomorphic to $M(n,n;\K)$, and hence $\cal O$ is isomorphic
as a homogeneous bundle to $\Gl(n,\A)$. The remaining task
of calculating the Lts of the symmetric space
$\cal O$ becomes easier by transforming everything with
the ``Cayley transform'' $R_n$: 
the Cayley transformed version of $\gl(n,\D)$ is
$$
R_n M(n,n;\D) R_n^{-1} = \Big\{ \pmatrix{a & \eps b \cr \eps c & d \cr} | \,
a,b,c,d \in M(n,n;\K) \Big\}.
$$
The Lie algebra of $\Gl(n,\A)$ is imbedded as the subalgebra given by the 
conditions $a=d$ and $b=c$.
The tangent spaces of our
two special orbits inside $R_n {\cal O} R_n^{-1}$ are complementary to this
subalgebra, namely
$$
\m = \Big\{ \pmatrix{Y & 0 \cr 0 & -Y \cr} | \, Y \in M(n,n;\K) \Big\}, \quad
V = \Big\{ \pmatrix{0 & \eps X \cr -\eps X & 0 \cr} | \, X \in M(n,n;\K) \Big\},
$$
so that $\m \oplus V$ is the Lts of $\cal O$, where the triple product is the
usual triple Lie bracket of matrices since the group action is by ordinary
conjugation of matrices. The formula shows that $\m$ clearly is isomorphic to
the usual Lts of $\Gl(n,\K)$, that $V$ is abelian and that $\Gl(n,\K)$ acts on 
$V$ by conjugation, i.e., the $r$-component of the corresponding Lts-representation of
$\m$ is the usual one (corresponding to the fact that $\cal O$ is isomorphic to 
$\Gl(n,\K[\eps])$ as a homogeneous bundle). However,
the whole
Lts representation just defined is not equivalent to the adjoint Lts representation of
$\gl(n,\K)$ on $\eps \gl(n,\K)$.
In fact, the corresponding involutive Lie algebras
$(\gl(n,\D),\gl(n,\A))$ and $(\gl(n,\A) \times \gl(n,\A), dia)$ are not isomorphic:
already for $n=1$, they are not isomorphic since 
$\gl(1,\A) \times \gl(1,\A) = \A \times \A$
is commutative, whereas $\gl(1,\D) = \D$ is not. 
\qed

\nin As mentioned above, instead of dual numbers
we could have taken for $\A$ another ring extension of the form
$\A_\mu = \K[X]/(X^2 - \mu)$ with arbitrary $\mu \in \K$ instead of $\mu=0$.
Then the symmetric space $\Gl(n,\K)$ has two different scalar extensions from
$\K$ to $\A_\mu$: the ``straight one'', simply gotten by taking the group case
$\Gl(n,\A_\mu)$, and another, ``twisted'' one, given by the homogeneous symmetric space
 $\Gl(n,\D_\mu)/\Gl(n,\A_\mu)$, where $\D_\mu$ is
the split Cayley-Dickson extension of $\A_\mu$ (see [Did06] for details). 
One could even replace here the algebra of square matrices by any other
associative $\K$-algebra.
For $\mu=-1$, we are back
in the example of the ``twisted complexification of $\Gl(n,\R)$''.
It is interesting that the interpretation of $\cal O$ as the 
``space of complex structures'' (i.e., endomorphisms with $E^2 = \mu$)
works only for invertible scalars $\mu$, whereas the interpretation given 
here works uniformely for all scalars.

\msk
\nin
{\bf 6.3. Jordan-extensions.} Besides general linear groups,
for all other ``classical'' symmetric spaces, there exist similar 
descriptions of symmetric bundle structures on the tangent bundle, see
[Did06] for the case of Grassmannians, Lagrangians and orthogonal groups.
The latter example gives rise to the ``$\D$-unitary groups'', analogues
of $\Sp(n)$ with $\H$ replaced by $\D$. 
%
%(By the way, this makes
%appear $\Sp(n,\K)$ as a unitary subgroup of $\Gl(n,M(2,2;\K))$...) 
%
The general construction principle behind these examples uses Jordan theory:

\Definition 6.4.
A {\it Jordan-extension} of a Lts $(\m,R)$ is given by a Jordan triple product
$T:\m^3 \to \m$ such that
$$
R(X,Y)Z=T(X,Y,Z)-T(Y,X,Z).
$$
Recall that a Jordan triple system (Jts) is a linear space $\m$ together with
a trilinear map $T:\m^3 \to \m$ such that

\ssk
\item{(JT1)} $T$ is symmetric in the outer variables: $T(u,v,w)=T(w,v,u)$
\item{(JT2)} 
$T(u,v,T(x,y,z))=T(T(u,v,x),y,z) -T(x,T(v,u,y),z) + T(x,y,T(u,v,z))$

\ssk
\nin For any Jts $T$, the trilinear map defined by
 $R_T(x,y)z=T(x,y,z)-T(y,x,z)$ is a Lts: we call the correspondence $T \mapsto
R_T$ the {\it Jordan-Lie functor} (cf.\ [Be00]).

\Theorem 6.5.
{\rm (The twisted regular representation defined by a Jordan extension.)}
Assume $(\m,R)$ is a Lts having a Jordan extension $T$.
Let $(\eps \m, r,m)$ be the regular representation of $\m$ (recall that
the corresponding split null extension $\m \oplus \eps \m$ is just the scalar
extension by dual numbers).
Then there exists another representation $(\eps \m,\tilde r,\tilde m)$
of $\m$,  in general not isomorphic to  the regular representation, but
  having the same $r$-component  (i.e., $\tilde r=r$). 

\Proof. 
We follow the lines of the proof of the corresponding statement for complexifications
in [Be00, Ch.\ III]:
let $\m[\eps] = \m \oplus \eps \m$ and $T[\eps]: \m_\eps^3 \to \m_\eps$
be the $\epsilon$-trilinear scalar extension of $T$ by dual numbers.
Then the {\it conjugation}
$$
\tau:\m[\eps] \to \m[\eps], \quad \tau(x+\epsilon y) =
\overline{x+\epsilon y} = x-\epsilon y 
$$
is a $\K$-automorphism of $T[\epsilon]$.
But for any involutive $\K$-automorphism, the ``$\tau$-isotope''
$$
\tilde T(u,v,w):=T[\epsilon](u,\tau(v),w)
$$
is again a Jts (cf.\ [Be00, Lemma III.4.5] for the easy proof). Moreover,
this new Jts is $\eps$-linear in the outer variables and
$\eps$-antilinear in the inner variable, and
since $\tau$ acts trivially on $\m$, restriction of this new Jts to
$\m^3$ gives us back $T$ again. Now we let
$$
\tilde R(X,Y):=R_{\tilde T}(X,Y)= \tilde T(Y,X) - \tilde T(X,Y).
$$
This is a Lts which coincides with $R$ on $\m$ since $T$ was chosen
to be a Jordan extension of $R$.
Next, $\epsilon \m$ is an ideal of $\tilde R$:
by (anti-)linearity it is an ideal of $\tilde T$, and hence it is one
of $\tilde R$.
Finally, if two terms belong to $\epsilon \m$, then
application of $\tilde T$ gives zero, and therefore also
application of $\tilde R$ gives zero.
Thus $\tilde \m$ is a Lts having the properties from Lemma 2.2, and hence
is the split null extension corresponding to a representation $(\tilde r,
\tilde m)$ on $\eps \m$.

\ssk
Now we show that $r=\tilde r$: for $x,y,v \in \m$,
$$
\tilde r(x,y,\eps v)=\tilde R(x,y,\eps v)=\tilde T(x,y,\eps v)-
\tilde T(y,x,\eps v)=r(x,y,\eps v).
$$
In order to prove that $(r,m)$ and $(\tilde r,\tilde m)$ are in general not
isomorphic, observe that the split null extension of $(r,m)$, being just
scalar extension by $\K[\eps]$, has the property that $R(\eps X,Y)=R(X,\eps Y)$.
On the other hand, 
$$
\tilde R(\eps X,Y)=\tilde T(\eps X,Y,\cdot)-\tilde T(Y,\eps X,\cdot)=
\eps (\tilde T(X,Y) + \tilde T(Y,X)) =
- \tilde R(X,\eps Y).
$$
Thus $R$ and $\tilde R$ cannot belong to isomorphic representations
 unless they vanish both. (This does not exclude that, as Lie triple systems over $\K$,
they may be isomorphic in special cases.)
\qed

\msk \nin {\bf 6.6. Final comments.}
Essentially, all classical Lie triple systems (and about half of the exceptional
ones) admit Jordan extensions
(cf.\ [Be00, Chapter XII]);
for instance, it is easily checked that $M(n,n;\K)$ with the triple product
 $T(u,v,w)=uvw+wvu$ is a Jts, and then
$$
T(u,v,w)-T(v,u,w) = uvw + wvu - (vuw + wuv) = [u,v]w - w[u,v] = 
 [[u,v],w],
$$
so that we have a Jordan extension of $\gl(n,\K)$. 
Correspondingly, $\Gl(n,\K)$ and essentially all classical symmetric spaces
admit on their tangent bundle a symmetric bundle structure that is different
from the usual one. 
We conjecture that (at least for simple finite-dimensional Lts over $\R$ or $\C$)
all symmetric bundle structures on the tangent bundle are exactly of two types

\ssk
\item{(1)} ``straight'': given by the canonical symmetric structure on $TM$, corresponding
to the regular representation of the Lts $\m$,
\item{(2)} ``twisted'': given by the construction from the  preceding theorem.

\ssk \nin
This conjecture is of course supported by the corresponding fact for 
{\it complexifications} of symmetric spaces, which (for simple finite-dimensional
Lts over $\R$) are either straight or twisted ([Be00, Cor.\ V.1.12]).
However, the proof given in loc.\ cit.\ for the complex and para-complex cases
 does not carry over to the tangent case (the invertibility of $i$, resp.\ $j$,
in $\K[i]$, resp.\ $\K[j]$
is used at a crucial point, and $\eps$ clearly is not invertible in $\K[\eps]$).
A proof covering all three cases at the same
time would be of great value for a better understanding of the ``Jordan-Lie functor''
(see [Be00]), and it should relate the ``extension problem for the Jordan-Lie functor''
with the extension problem for Lts representations as discussed here. 
One might conjecture that an interpretation in terms of the Cayley-Dickson process,
which turned out to be useful in the special case of $\Gl(n,\K)$, could be the key
for proving the conjecture, but this is not clear at present.

\bigskip
\nin {\bf Acknowledgement.}
The first named author would like to thank the Hausdorff Institute (Bonn) for
hospitality when part of this work was carried out.

\def\entries{

% \[Ar66 Artin, E., {\it Geometric Algebra}, Interscience, New York 1966

%\[B57 Berger, M., ``Les espaces sym\'etriques non compacts", Ann. Ec. 
%Norm. Sup. (3) 74 (1957),  85--177

\[Barr96
Barr, M., ``Cartan-Eilenberg cohomology and triples'',
J. Pure Applied Algebra {\bf 112} (1996), 219 -- 238

\[Beck67 Beck, J.\ M., {\it Triples, algebra and cohomology}, thesis,
Columbia University 1967 (available at
http://www.emis.de/journals/TAC/reprints/index.html )

\[Be00
Bertram, W., {\it The geometry of Jordan- and Lie structures},
Springer Lecture Notes in Mathematics {\bf 1754},  Berlin 2000

\[Be02
 Bertram, W., ``Generalized projective geometries: general
theory and equivalence with Jordan structures", Adv. Geom. {\bf 2}
(2002), 329--369 
%(electronic version: preprint 90 at 
%http://www.uibk.ac.at/mathematik/loos/jordan/index.html)

\[Be08
Bertram, W., {\it Differential Geometry, 
Lie Groups and Symmetric Spaces over General Base Fields and Rings.}
Mem.\   AMS {\bf 900}, 2008;  math.DG/ 0502168 

%\[Bou Bourbaki, N., {\it Groupes et alg\`ebres de Lie},

%\[Cl92 Clerc, J.-L., ``Repr\'esentations d'une alg\`ebre de Jordan,
%polyn\^omes invariants et harmoniques de Stiefel'',
%J. reine angew. Math. {\bf 423} (1992), 47--71

\[Did06 Didry, M., {\it Structures alg\'ebriques associ\'ees aux espaces
sym\'etriques}, thesis, Institut Elie Cartan, Nancy 2006
(see http://www.iecn.u-nancy.fr/ $\tilde{ }$ didrym/ )

\[Did07 Didry, M., ``Construction of groups associated to Lie- and
Leibniz algebras'', J.\ Lie theory {\bf 17}, no.\ 2 (2007), 399 -- 426 

\[Ei48  Eilenberg, S., ``Extensions of general algebras'',
Ann. Soc. Pol. Math. {\bf 21} (1948),  125--134

\[FK94 Faraut, J. et A. Kor\`anyi, {\it Analysis on Symmetric Cones},
Clarendon Press, Oxford 1994

\[Ha61 Harris, B., ``Cohomology of Lie triple systems and Lie algebras
with involution'', Trans. A.M.S. {\bf 98} (1961),  148 -- 162

\[Ho02
Hodge, T.L. and B.J. Passhall, 
``On the representation theory of Lie triple systems.''
 Trans. Amer. Math. Soc. {\bf 354} 2002 N. 11 ,
 4359--4391

\[Jac51 Jacobson, N., ``General representation theory of Jordan algebras'',
Trans. A.M.S {\bf 70} (1951), 509--530

\[KoNo69
Kobayashi, S.\ and K.\ Nomizu, {\it Foundations of Differential Geometry. Volume
II.}, Wiley, New York 1969

%\[KT04
%Kubo, F. and Yo. Taniguchi. ``A controlling cohomology of the deformation theory
%of Lie triple systems.''
% J. Algebra, {\bf 278 (1)} (2004)  242 -- 250

\[Li52 Lister, W.G., ``A structure theory for Lie triple systems'',
Trans. A.M.S. {\bf 72}(1952), 217 -- 242

\[Lo67 Loos, O., ``Spiegelungsr\"aume und homogene symmetrische
R\"aume''
Math. Z.{\bf 99} (1967), 141 -- 170

\[Lo69 Loos, O., {\it Symmetric Spaces I}, Benjamin, New York 1969

\[Lo73 Loos, O., ``Representations of Jordan triples'', Trans. A.M.S.
{\bf 185} (1973), 199--211

\[Lo75 Loos, O., {\it Jordan Pairs}, Springer LNM {\bf 460}, Berlin 1975

%\[Mc71 McCrimmon, K., ``Representations of quadratic Jordan algebras'',
%Trans. A.M.S. {\bf 153} (1971), 279--305

%\[McC04
%Kevin McCrimmon. A taste of Jordan algebras. Universitext. Springer-Verlag, New York,
%2004.

\[MoPe06
Mostovoy, J. and  J.M.\ Pérez-Izquierdo,
``Ideals in non-associative universal enveloping algebras of Lie triple systems.''
arXiv:math/0506179

%\[Mo01
% Shigeyuki Morita. 
%Geometry of differential forms. American Mathematical Society, Providence,
%RI, 2001. Translated from the two-volume Japanese original (1997, 1998) by Teruko
%Nagase and Katsumi Nomizu, Iwanami Series in Modern Mathematics.

\[Sa99
Sabinin, Lev V., 
{\it Smooth quasigroups and loops}, 
Mathematics and its Applications, {\bf 492}. Kluwer Academic Publishers, Dordrecht, 1999

%\[Se65
%J-P. Serre. {\it Lie algebras and Lie groups},
%W. A. Benjamin, Inc., New York-Amsterdam, 1965

\[Sm05
Smirnov, O.N., ``Imbedding of Lie triple systems into Lie algebras'',
preprint, 2005

%\[Ya60
%Kiyosi Yamaguti. On the cohomology space of Lie triple system. Kumamoto J. Sci. Ser.
%A, 5 :44--52 (1960), 1960

%\[ZSSS82
%Zhevlakov,K.A., A. M. Slinko, I. P. Shestakov, and A. I. Shirshov. 
%{\it Rings that are nearly
%associative}, volume 104 of Pure and Applied Mathematics. Academic Press Inc. [Harcourt
%Brace Jovanovich Publishers], New York, 1982. 

}

\references

\lastpage

\vfill\eject

\end

{\bf Remarks and comments for future work:}

\bigskip
{\bf 1. Several symmetric bundle structures on $TM$ and associated tensors.}
From the point of view of connections, the uniqueness question now takes the
following form: assume $F_1$ and $F_2$ are two symmetric bundles over $M$,
so that, as reflection spaces (or as homogeneous bundles, if $M=G/H$), they 
agree with each other; in particular, as abstract bundles they are the same bundle
$F$ over $M$. By applying the tangent functor $T$, the tangent bundle $TF=TF_1=TF_2$ 
also carries two symmetric space structures. These symmetric spaces are uniquely
characterized by their canonical connections, say $L_1$ and $L_2$
 (which are linear connections on $TTF$ over $TF$).  
The difference $A:=L_1 - L_2$ is an algebra-field (tensor field of type $(2,1)$)
on $TF$. This field is a field of {\it commutative} algebras (since the canonical connection
is torsionfree), and it is
 invariant under the group $G=G(M)$ which preserves both
structures; hence it is enough to look at $A$ on the fiber $V=F_o$ over a base point
$o \in M$, i.e.\ to look at
$A_u:T_u (TF) \times T_u(TF) \to T_u (TF)$ for $u \in V$.
Since the splitting $T_u(TF)=V_u \oplus H_u \cong V \oplus \m$ is the same for
$F_1$ and $F_2$, we end up with a map
$$
(V \oplus \m) \times V \times  (V \oplus \m) \to  V \oplus \m, \quad
(X,u,Y) \mapsto A_u(X,Y).
$$
SHOW: this depends linearly on $u$, and the only non-zero component is a trilinear
map $\m \times V \times \m \to V$. Its properties should be rather close to those
defining a Jts (note that it must be symmetric in $X$ and $Y$).
The relation with the $m$-components of the two underlying Lts representations
of $\m$ is given by the Maurer-Cartan equation
$$
R_1(X,Y) = R_2(X,Y) + \nabla A((X,Y)-(Y,X)) + [A_X,A_Y],
$$
where $R_1,R_2$ are the curvature of $L_1,L_2$ and $\nabla$ is covariant
derivative with respect to $L_1$.
Possibly, an algebraic approach to these questions would involve some
cohomology theory for Lie triple systems, cf.\ [Ha61], [Ya60].

\msk
{\bf Conjecture.}
 We  conjecture that {\it symmetric bundle structures
on $TM$ other than the usual one are in one-to-one correspondence with Jordan
extensions $T$ of the curvature $R$ of $M$}.
In other words, we conjecture that all symmetric bundle structures on $TM$,
other than the usual one, are given by the construction from Proposition 5.2.

Some more evidence for this: first step: in the situation of Prop. 5.2., try
to recover the Jordan tensor in a geometric way (for instance, relate it to the
difference of the two different affine connections on $T(TM)$ over $TM$,
which is a $G$-invariant
commutative algebra field vanishing on the zero section and certainly depending
linearly on the fiber variable: it must be the Jordan tensor),
second step: try to reverse this heuristic and to define the Jordan tensor
in a geometric manner. It may happen that the old straight/twisted opposition
reappears and that the conjecture only holds when there are no other 
``complexifications'' -- here is the old ansatz for that:
Assume that $\q$ is simple. Then we would like to show that
$TM$ is either straight or twisted. Here the following could
be helpful: it straight iff
$$
\q \otimes \q \to \Hom(\q,\epsilon \q), \quad
(X,v) \mapsto R(X,\epsilon v) 
$$
is skew-symmetric, and twisted iff this map
is symmetric. It seems reasonable that in the 
simple case only these two possibilities survive...

10/07: try to follow the remark from the $\Gl(n)$-case: the direct interpretation
as ``space of complex structures'' does not carry over as soon as the scalar is
non-invertible, but the more subtle point of view of ``twisted $\K$-form of $\DD$''
does: can one adapt this to prove the conjecture??

Let us just add some calculations that point into this direction: the map $\tilde m$
from the preceding proof is given by, if $x,v,y \in \m$,
$$
\eqalign{
\tilde m(x,\epsilon v,y)=\tilde R(x,\epsilon v)y=&=
\tilde T(\epsilon v,x,y) - \tilde T(x,\epsilon v,y) \cr
& = \epsilon (T(v,x,y)+T(x,v,y)) \cr
& = \epsilon (2T(x,v,y) - T(x,v,y) + T(v,x,y)) \cr
& = \epsilon (2T(x,v,y) +  R(x,v,y)). \cr}
$$
Thus, with $m(x,y)v=\epsilon R(x,v)y$ belonging to the
straight structure, we have
$$
\tilde m(x,y) - m(x,y) = 2 \epsilon T(x,\cdot,y).
$$
Hence the Lie algebra of the stabilizer, acting on $ßm \oplus \eps \m$, is now
$$
\pmatrix{R(x,y) & \epsilon (T(a,b)+T(b,a)) \cr
0 & \epsilon R(x,y)\cr}, \quad
x,y,a,b \in \q,
$$
instead of the  $\epsilon r(a,b)$ in the upper right corner appearing in the
``straight'' case.
Note finally that the ``twisted'' symmetric
structure on $TM$ induces, via 4.2, also a new, ``twisted'' symmetric structure
on the cotangent bundle $T^*M$ and on all other bundles obtained by linear
algebra constructions.

\bigskip
{\bf Comments on the proof for $\Gl(n,\K)$ (Theorem 6.2):}
The proof gives, moreover, the following explicit formula for the
interesting $\m \times V \times \m \to V$-component of the Lie triple product
on $\m \oplus V$:
$$
\eqalign{
\Big[ \big[
\pmatrix{Y & 0 \cr 0 & -Y \cr}, & \pmatrix{0 & \eps v \cr -\eps v & 0 \cr} \big],
\pmatrix{Z & 0 \cr 0 & -Z \cr} \Big] = \cr
&
\eps \pmatrix{0 & YvZ+ZvY+vYZ+ZYv \cr
-(YvZ+ZvY+vYZ+ZYv) & 0 \cr}
\cr}
$$
We will see below that this formula has a Jordan-theoretic interpretation.

\bigskip
\nin {\bf 
Here old idea to work with a ``structure tensor'' as in [Be00]:
The $\epsilon$-tensor.}
Let $M$ be a manifold, $TM$ its tangent bundle and
$\pi:TM \to M$ the canonical projection.
For $v \in T_pM$ choose a splitting 
$T_v(TM)=H_v \oplus V_v$
of the tangent space
$T_v(TM)$ such that $V_v = \ker(T_v \pi) \cong T_pM$;
then
$$
T_v\pi: H_v \to T_pM
$$
is a linear bijection. 
Define
$$
\epsilon_v: T_v(TM)=H_v \oplus V_v \to T_v(TM), \quad
(h,x) \mapsto (0,T_v(\pi)h).
\eqno (4.1)
$$
The map $\epsilon_v$ is independent of the complement $H_v$
of $V_v$, and thus $\epsilon = (\epsilon_v)_{v \in TM}$
is a well-defined tensor field of type $(1,1)$ on $TM$, called
the {\it $\epsilon$-tensor}.
In matrix form we have, identifying 
$T_v M = H_v \oplus V_v$ with $T_p M \oplus T_p M$, 
$$
\epsilon_V =
\pmatrix{0 & 0 \cr \1 & 0 \cr}.
\eqno (4.2)
$$
Clearly, $\epsilon^2 =0$, and $\ker \epsilon_v = \ker( T_v\pi)$.
The $\epsilon$-tensor can be seen as a dual number analog
of the almost complex structure of a complex manifold.
{\bf Straight and twisted $\epsilon$-tensors on symmetric
bundles.}
Now assume
$M = G/H$ is a symmetric space and consider a structure of
symmetric bundle $TM = L/B$ on its tangent bundle (with
$B$ assumed to be connected).
Then the $\epsilon$-tensor is $L$-invariant: in fact,
$B$ is
generated by the exponentials of the matrices 
$R(X+v,Y+w)$ given by Eqn. (2.4) (note that here
$r(X,Y)=R(X,Y)$); but it is immediately verified
that the matrices (4.2) and (2.4) commute, that is, the condition
$$
R(X,Y) \epsilon Z = \epsilon R(X,Y)Z
\eqno (4.3)
$$
holds for all vector fields $X,Y,Z$ on $TM$, where $R$ is the curvature
of $TM=L/B$.
We say that $(R,\epsilon)$ is 
{\it straight} if in addition the condition
$$
R(\epsilon X,Y) = R(X,\epsilon Y)
\eqno (4.4)
$$
holds,
and {\it twisted} if in addition the condition
$$
R(\epsilon X,Y) = -R(X,\epsilon Y)
\eqno (4.5)
$$
holds.
Using the Jacobi identity, it easily seen that
(4.3) and (4.4) together imply that $R$ is
$\epsilon$-trilinear (cf. [Be00, Lemma III.1.7]).
Thus this case corresponds to the canonical symmetric
structure on the tangent bundle defined in Section 2.10.
On the other hand, (4.3) and (4.5) together imply
that $\epsilon $ is a derivation of $R$.
Similarly as in [Be00, Ch. III], we define the
associated {\it structure tensor} by
$$
T(X,Y,Z):={1\over 2}(\epsilon R(X,Y)Z-R(X,\epsilon Y)Z).
\eqno (4.6)
$$
Since $\epsilon$ and $R$ are invariant tensor fields on $TM$,
so is $T$. The proofs of the following properties are
similar to those for complex structures ([Be00, Prop. III.2.3.
Prop. III.2.4]) and are therefore omitted here:

\ssk

\item{(1)}
$T$ is symmetric in the outer variables:
$T(X,Y,Z)=T(Z,Y,X)$,
\item{(2)}
$T$ is related to $R$ via
$T(X,Y,Z)-T(Y,X,Z) = -\epsilon R(X,Y)Z$,
\item{(4)}
$T$ satisfies the identity
(where we let $T(X,Y)Z:=T(X,Y,Z)$)
$$
T(X,Y)T(U,V,W)=
T(T(X,Y)U,V,W)-T(U,T(Y,X)V,W) + T(U,V,T(X,Y)W).
\eqno {\rm (JT2)}
$$

\nin
(1) and (4) say that $T$ is a {\it Jordan triple system}.
Since $\epsilon \q$ is an ideal, $T$ takes values
in $\epsilon \q$ {\bf and thus all terms in (JT2) are
zero....not very interesting !  need that
$T: \q \times \q \times \q \to \q \cong \epsilon \q$
is a JTS: not the same thing !}
Thus define, with $V=\q$, for $X,Y,Z \in \q$,
$$
T(X,Y,Z)=
{1\over 2}(\epsilon R(X,Y)Z-R(X,\epsilon Y)Z) =
{1 \over 2}(r(X,Y)Z-m(X,Z)Y).
$$
Then $T(X,Y,Z)=T(Z,Y,X)$ remains valid.
We have
$$
T(X,Y)-T(Y,X)=r(X,Y) = \epsilon R(X,Y)=R(X,Y)
$$
since $m(X,\cdot)Y$ is symmetric in $X,Y$.
This remains a derivation of everything.
Next,
$$
T(X,Y)+T(Y,X)=-m(X,\cdot)Y=-R(X,\epsilon \cdot)Y.
$$
Here is a problem: this map goes originally from
$\q$ to $\epsilon \q$. I don't know how to write
it as a map $\q \to \q$...
If we could show that it is an antiderivation of $T$
we would have (JT2) !
-- Other idea:
if $\epsilon$ is twisted, then $1+\epsilon$ is an
automorphism which is an infinitesimal automorphism (!),
 it defines thus  a section of the
tangent bundle; call it the {\it Euler operator} $E$.
Aim: verify the properties from [Be00] of an 
``Euler operator on $M$'', that is,
$\q^b = \q + [E,\q]$ is essentially a twisted polarized Lts.

\msk
{\bf General extension problem: algebraic approach.} [Cf. Harris, McCrimmon,...]
Given $r$, we want to classify all possible extensions
$m$ up to isomorphism.
A necessary condition for the existence of $m$ is that
the given sequence $0 \to V \to \tilde \q \to \q \to 0$
splits algebraically.
Next condition:
the map $r$ has to be in the image of the antisymetrization
map
$$
\Hom(\q \otimes \q,\End(V))^\h \to \Hom(\Lambda^2 \q,\End(V))^\h.
$$
Next condition:
the right-hand side of (R4) is antisymmetric in $X$ and $W$.
That is, for all $D \in \h$, the map
$$
\q \otimes \q \to \End(V), \quad
X \otimes W \mapsto m(X,D \cdot W) - \rho(D) m(X,W)
$$
is antisymmetric.
The space of $m$'s satisfying this condition is a linear
subspace $U$ and thus is never empty.
If two $m$'s are given, the difference satisfies again this
condition. If the space $U$ is trivial, then there
is at most one extension.
Intersecting $U$ with the preimage of $r$ under antisym.
yields an affine subspace which already could be rather
small...
However, we have not yet fully expressed the condition
(R4). The problem is that it is quadratic in $m$.

Assume $m,m':\q \otimes V \otimes \q \to V$ are two
symmetric space extensions of the same $\h$-module $\q$.
(Here probably already we have to be more careful:
the complement $\q$ of $V$ in $\tilde \q$ depends on the
symmetric space structure, therefore we cannot take the
same $\q$...)
Then
$$
f:=m-m':\q \otimes \q \to \End(V)
$$
is symmetric by (R2).
(R3) means that it is $\h$-equivariant.
(R4) implies for $f$ that... 
(more difficult: dependence on $f$ is quadratic on the right hand
side and linear on the left hand side.)
Condition on $f$ to be a cocylce ??
kernel of some differential.
problem: space of $f$'s is not a linear nor an affine space...

Next:
coboundaries: condition on $f$ corresponding to $m$ and $m'$
being isomorphic?

\msk \nin
{\bf The universal bundle.}
%
%[NOT YET CLEAR!]
Let $(\m,R)$ be a Lie triple over $\K$.
A {\it universal representation of $\m$} is a module
$(r,m)$ in an associative algebra $A$ with identity
$e$  such that every other representation
$(r',m')$
factors via a homomorphism $\phi:A \to A'$ of
{\it associative algebras}:
$$
r'=\phi \circ r:\m \times \m \to A', \quad
m'=\phi \circ m:\m \times \m \to A'.
$$
Such an object is unique and exists  (cf. [Lo73], [Lo75] with reference to
Cohn, ``Univ. algebra''), denote it by ${\bf u}(\m)$.
Note that $(r,m)$ followed by the left regular representation
is a representation on $A$ in the previous sense.
In any case, there is now an (infinite dimensional!)
symmetric bundle $\bf A$ over $M$ with typical fiber $A$ and which is universal
in the sense that
for any symmetric bundle $F$ over $M$ there is a homomorphism
$\Phi:{\bf A} \to F$.

 Essentially, $A = U^-:=\End({\bf u}(\g)^-)$ is the universal algebra
considered in [Ha61], associated to the standard imbedding $\g$. Since the standard
imbedding is not functorial, in general we just get a hom of ${\bf u}(\m)$ into
this object, but not in the other direction. However, this hom should not be ``too far''
from being an isomorphism. 

Note: $U^-$ is infinite-dimensional: take PBW-basis, those with an odd number of
tensors from $\g^-$ form a basis of $U^-$, and there is an infinite number of them.

Problem: what is the relation with the ``non-associative
enveloping algebra'' from [MoPe06] ?

\msk
\nin {\bf  The cotangent bundle.}
Since $TM$ is a symmetric bundle, 3.2 implies that the
cotangent bundle $T^*M$ also is a symmetric bundle.
Using 3.5 or 3.6 we can also define the analog of 
the bundle $T^*M \otimes TM$: it is the bundle with
typical fiber
$$
\Hom(\h,\m) \oplus \Hom(\m,\h)
$$
which comes from the $(\g,\sigma)$-module
$\End(\g)=\gl(\g)$ with involution $\sigma_*$.
Additional question: is $T^*M$ a symplectic symmetric
space ? No, according to Lionel (it is pseudo-Riemannian).
Try to explain.
General def of invariant forms for modules: must
$\q$ and $V$ be always isotropic (unless the sum is direct) ?
this would reduce possibitlities to $\dim \q = \dim V$...

\sectionheadline{X. Representations of Jordan pairs and triples}

The theory of general representations of Jordan pairs and
triples is
 due to Loos [Lo73], [Lo75].
We briefly present the linear version:
we write the JTS $T$ on $\q$ as 
$$
T(X,Y,Z)=T(X,Y)Z=P(X,Z)Y
$$
with the usual identities (JT1), (JT2).
Then for $\q \oplus V$ with
$$
t:\q \times \q \to \End(V), \quad
p:\q \times \q \to \End(V)
$$
 to be  a split
null extension we need the identities

\item{(1)}
$p(X,Y)=p(Y,X)$
\item{(2)}
$[t(X,Y),t(U,V)]=t(T(X,Y)U,V) - t(U,T(Y,X)V)$,
\item{(3)}
$p(T(X,Y)\cdot U \otimes V)=
t(X,Y)p(U,V)+p(U,V)t(Y,X)$,
\item{(4)}
$p(X,T(U,V)W)=
t(W,V)p(X,U)-p(U,W)t((V,X)+t(U,V)p(X,W)$

\ssk
\nin
In good cases (semisimple e.g.), (2) means that $t$ passes
as a homomorphism
$$
\str(T) \to \gl(V), \quad T(X,Y) \mapsto t(X,Y)
$$
of Lie algebras.
Thus $t$ simply corresponds to homogeneous bundles 
over the twisted para-complexification
of $M$ (and by integrability also over the compactification of
$M$).
A compatible choice of $p$ equips this bundle
with the structure of a symmetric space with twist
(i.e. the split null extension is again a JTS).

Problems:
-- Interprete all this on the level of generalized
projective geometries over a general base ring $\K$;

-- express the Jordan triple product in terms of
$T(TM)$ with respect to $\epsilon_1$ twisted and
$\epsilon_2$ straight.

--
Jordan algebras are more complicated objects; relate
things to null-systems:

\sectionheadline{XX. Representations of Jordan algebras}

\nin {\bf Special and general Jordan algebras.}
A {\it special Jordan algebra} is a subspace of an associative
algebra (over a field $\K$ which we assume to  be of characteristic
different from 2 or 3) which is closed under the composition 
$x \mapsto x^2$; it is a subalgebra with respect to the product
$$
x  y = {x\cdot y + y\cdot x \over 2}.
\eqno (3.1)
$$
Let $L_a(x):= a  x=xa$. Then in any special Jordan algebra the
following identites are satisfied:

\ssk
\item{(1)}
$[[L_a,L_b],L_c] = L_{a(bc)-b(ca)}$,
\item{(2)}
$[L_{ab},L_c] + [L_{bc},L_a] + [L_{ca},L_b]=0.$

\ssk
\nin
If we let $a=b=c=:x$, then (1) becomes trivial and (2) becomes
$$
[L_{x^2},L_x] = 0. \eqno (J2)
$$
One can show that (J2) in turn implies (1) and (2) which are 
in fact the complete linearization of (J2).
A {\it (general) Jordan algebra} over $\K$ is an algebra with
a commutative product satisfying (J2).

\msk
\nin {\bf Special representations.}
If $V$ is a Jordan algebra, then a {\it special representation}
is a homomorphism
$$
\phi: V \to \End(E)
$$
into an endomorphism algebra with the product (3.1).
If $V$ has a unit $e$, then one requires also that $\phi(e)=\id_V$.
We define on
$$
\tilde V := V \oplus E
$$
a product by
$$
(v \oplus a)(w \oplus b):=vw \oplus (\phi(v)b + \phi(w)a).
$$
Then, for $v \in V$, $a \in E$, 
$$
\eqalign{
(v,0)^2 ((v,0)(0,b))&
=(v^2,0)(0,\phi(v)b)=(0,\phi(v^2)\phi(v)b)
\cr
& = (0,\phi(v)^3b)= (v,o)((v,0)^2 (0,b)).
\cr}
$$
This proves (J2) in $\tilde V$ since (J2) is trivial in all other cases.

\msk
\nin {\bf General representations.}
It is not possible to reverse the preceding calculations:
if we assume that $\tilde V = V \oplus E$ is a Jordan algebra
with an ideal $E$ on which the product is zero, then we cannot
conclude that
$$
\phi: V \to \End(E), \quad v \mapsto (x \mapsto vx)
$$
is a special Jordan algebra representation (the Jordan identity
$x^2(xv)=x(x^2 v)$ is of degree 3 and is not sufficient to
deduce the second order identity $\phi(x^2)=\phi(x)^2$,
i.e., $x^2 v=x(xv)$ which we would like to have).
Therefore we {\it define} a general representation of a Jordan
algebra $V$ to be a map
$$
L: V \to \End(E), \quad v \to L(v)
$$
such that $[L(v)^2,L(v)]=0$ for all $v \in V$.
Thus the special representations (which are defined by the identity
$L(v^2)=L(v)^2$) are general ones, but the converse is false:
in fact, the {\it regular representation}
$$
V \to \End(V), \quad v \mapsto L(v):=L_v
$$
of a Jordan algebra on itself is in general not special
(it is special iff $v(vx)=(vv)x$ iff $v(wx)=(vw)x$ iff
it is associative and commutative). 

Proposition:
If $E,F$ are general representations of a Jordan algebra $V$,
then so are $E \oplus F$ and $E^*$ but not (!?)
 $\Hom(E,F)$, $E \otimes F$.

Question: Space interpretation of all this? On two levels:
generalized projective geometries with inner polarity;
quadratic prehomogeneous symmetric spaces.
For the latter case, (1) implies
$$
L([a,b,c])=[L(a),L(b),L(c)]
$$
(where $[a,b,c]=a(bc)-b(ac)$ is the associator in $V$ which is
a Lie triple on $V$) and we get a homomorphism
$$
M \to \Gl(E)
$$
of symmetric spaces. In fact, if $E=V$ is the regular rep.,
 this is the quadratic representation
$Q:M \to \Gl(V)$. 

\msk
{\it In Null.tex make also theory of prehomogeneous symmetric
spaces: $M \subset V$, aff. structure of $V$ invariant under
$G(M)$,
 $Q:V \to \End(V)$ quadratic,
$\mu(x,y)=P(x) P(y)^{-1}y$ symm. space, without unit.}

\msk
\nin
{\bf Brainstorm (old file): Questions:}

1. Duality principle: what is the global version of the 
duality principle from [Lo73] and [Lo75]?
Note  the striking similarity of (PG1), (PG2)
and identities of the dual representation...

2. Cotangent bundle: symplectic and pseudo-Riemannian form;
relation with twisted polarized spaces?

3. Universal envelopes:
Universal envelope in Lie theory corresponds to the bundle
$T^\infty M$; in Jordan theory this cannot be the case since
there universal envelopes are finite-dimensional...
 
In general, higher order tangent bundles (sense of Loos, Pohl!)
 are not symmetric
spaces in a natural way; is there nevertheless some significant
relation with differential operators? 
with the Cartan connection?

4. Sections: If $E$ is a symmetric space bundle over $M$, what
is the role of the homomorphisms $M \to E$ which are sections?
If $E=TM$, then these are the derivations.
For $E=T^*M$?

5. The bundle $S^2 TM$ and invariant pseudo-metrics.
(attn. $S^2TM$ is not symmetric...)

6. Duflo iso for symmetric spaces and gen. proj. geos...

7. The bundles whose sections are nominators, resp. denominators...

\msk
\nin {\bf Group cases and specializations.}
If $V$ is finite dimensional, the general linear group 
$\Gl(V)$ is a symmetric space with product map
$$
\mu(g,h)=g h^{-1} g
$$
and associated Lie triple $\gl(V)$ with
$$
[X,Y,Z] = {1 \over 4} [[X,Y],Z].
$$
Even if $V$ is of arbitrary dimension we will equip $\Gl(V)$ and
$\gl(V)$ with these products.
A {\it specialization of a symmetric space $M$}, resp. of a LTS
$\q$, is a homomorphism
$$
\pi : M \to \Gl(V), \quad {\rm respectively} \quad
\pi : \q \to \gl(V).
$$
Slightly more general:
a {\it group specialization of $M$} is a homomorphism into
an (abstract) group. Then every symmetric space has an
``adjoint specialization'', namely $x \mapsto s_x$
(or $x \mapsto s_x s_o = Q(x)$).
{\it
(Question: why is there no such thing for gen. proj. geos.?
Answer: we need some sort of central extension in order
to define a structure on $\Aut(X,X')$...)
}

\vfill \eject

\bye